\documentclass[12pt]{amsart}

\usepackage{epsf}
\usepackage{graphicx}
\usepackage{latexsym}
\usepackage{amsfonts}
\usepackage{amsmath}
\usepackage{amssymb}
\usepackage{amscd}
\usepackage{amsthm}
\usepackage{bm}

\def\R{\mathbb{R}}

\def\Z{\mathbb{Z}}
\def\Q{\mathbb{Q}}

\def\CFK {\widehat{\operatorname{CFK}}}
\def\HFK {\widehat{\operatorname{HFK}}}
\def\HF {\widehat{\operatorname{HF}}}
\def\D {\mathcal{D}}
\def\QQ {\mathcal{Q}}
\def\RR {\mathcal{R}}
\def\SS {\mathcal{S}}
\def\T {\mathbb{T}}
\def\s {\mathfrak{s}}
\def\srel {\underline{\mathfrak{s}}}
\def\x{\mathbf{x}}
\def\y{\mathbf{y}}
\def\z{\mathbf{z}}
\def\w{\mathbf{w}}

\newcommand{\gen}[1] {\left\langle #1 \right\rangle}
\def\Th{^{\text{th}}}

 \DeclareMathOperator{\lk}{lk}
\DeclareMathOperator{\Sym}{Sym} \DeclareMathOperator{\Spin}{Spin}
\DeclareMathOperator{\Span}{Span}

\def\minus{\smallsetminus}

\newtheorem{prop}{Proposition} [section]
\newtheorem{lemma}[prop]{Lemma}
\newtheorem{thm}{Theorem}

\newtheorem{conj}[prop]{Conjecture}

\theoremstyle{definition}

\theoremstyle{remark}
\newtheorem{remark}[prop]{Remark}

\begin{document}
\title
 [Computing Knot Floer Homology in Cyclic Branched Covers]
 {Computing Knot Floer Homology in Cyclic Branched Covers}

\author{Adam Simon Levine}
\address {Department of Mathematics, Columbia University \\ New York, NY 10027}
\email {alevine@math.columbia.edu}

\begin{abstract}
We use grid diagrams to give a combinatorial algorithm for computing
the knot Floer homology of the pullback of a knot $K$ in its
$m$-fold cyclic branched cover $\Sigma^m(K)$, and we give
computations when $m=2$ for over fifty three-bridge knots with up to
eleven crossings.
\end{abstract}

\maketitle

\section{Introduction}

Heegaard Floer knot homology, developed by Ozsv\'ath and Szab\'o
\cite{OSzKnot} and independently by Rasmussen \cite{R1}, associates
to a knot $K$ in a three-manifold $Y$ a bigraded group $\HFK(Y,K)$
that is an invariant of the knot type of $K$. If $K$ is a knot in
$S^3$, then the inverse image of $K$ in $\Sigma^m(K)$, the $m$-fold
cyclic branched cover of $S^3$ branched along $K$, is a
nulhomologous knot $\tilde K$ whose knot type depends only on the
knot type of $K$, so the group $\HFK(\Sigma^m(K), \tilde K)$ is a
knot invariant of $K$. In this paper, we describe an algorithm that
can compute $\HFK(\Sigma^m(K), \tilde K)$ (with coefficients in
$\Z/2$) for any knot $K \subset S^3$, and we give computations for a
large collection of knots with up to eleven crossings.

Any knot $K \subset S^3$ can be represented by means of a \emph{grid
diagram}, consisting of an $n \times n$ grid in which the centers of
certain squares are marked $X$ or $O$, such that each row and each
column contains exactly one $X$ and one $O$. To recover a knot
projection, draw an arc from the $X$ and the $O$ in each column and
from the $O$ to the $X$ in each row, making the vertical strand pass
over the horizontal strand at each crossing. We may view the diagram
as lying on a standardly embedded torus $T^2 \subset S^3$ by making
the standard edge identifications; the horizontal grid lines become
$\alpha$ circles and the vertical ones $\beta$ circles. Manolescu,
Ozsv\'ath, and Sarkar \cite{MOS} showed that such diagrams can be
used to compute $\HFK(S^3,K)$ combinatorially; we shall use them to
compute $\HFK(\Sigma^m(K), \tilde K)$ for any knot $K \subset S^3$.
(See also \cite{BG, MOSzT, SW}.)

Let $\tilde T$ be the surface obtained by gluing together $m$ copies
of $T$ (denoted $T_0, \dots, T_{m-1}$) along branch cuts connecting
the $X$ and the $O$ in each column. Specifically, in each column, if
the $X$ is above the $O$, then glue the left side of the branch cut
in $T_k$ to the right side of the same cut in $T_{k+1}$ (indices
modulo $m$); if the $O$ is above the $X$, then glue the left side of
the branch cut in $T_k$ to the right side of the same cut in
$T_{k-1}$. The obvious projection $\pi: \tilde T \to T$ is an
$m$-fold cyclic branched cover, branched around the marked points.
Each $\alpha$ and $\beta$ circle in $T$ intersects the branch cuts a
total of zero times algebraically and therefore has $m$ distinct
lifts to $T$, and each lift of each $\alpha$ circle intersects
exactly one lift of each $\beta$ circle. (We will describe these
intersections more explicitly in Section \ref{sec:grid}.)

\begin{figure} \label{fig:differential}
\includegraphics[scale=0.75]{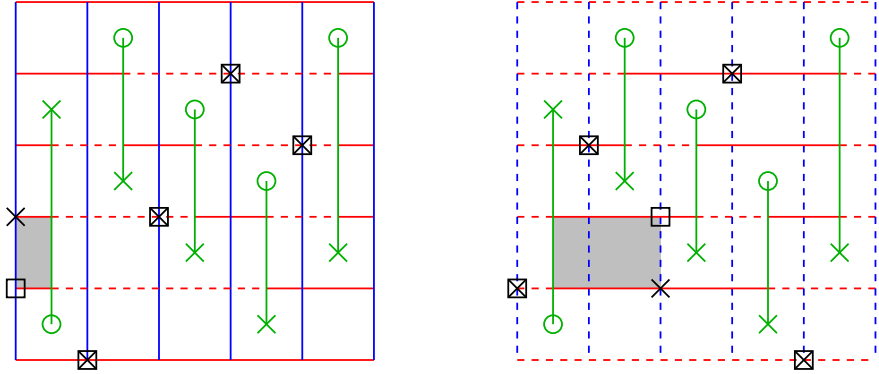}
\centering \caption{Heegaard diagram $\tilde D = (\tilde T, \tilde{
\bm\alpha}, \tilde{\bm\beta}, \tilde\w, \tilde\z)$ for
$(\Sigma^2(K), \tilde K)$, where $K$ is the right-handed trefoil.
The solid and dashed lines represent different lifts of the $\alpha$
(horizontal/red) and $\beta$ (vertical/blue) circles. The black
squares and crosses represent two generators of $\tilde C =
\CFK(\tilde\D)$, and the shaded region is a disk that contributes to
the differential.}
\end{figure}

Denote by $\RR$ the set of embedded rectangles in $T$ whose lower
and upper edges are arcs of $\alpha$ circles, whose left and right
edges are arcs of $\beta$ circles, and which do not contain any
marked points in their interior. Each rectangle in $\RR$ has $m$
distinct lifts to $\tilde T$ (possibly passing through the branch
cuts as in Figure \ref{fig:differential}); denote the set of such
lifts by $\tilde\RR$.

Let $\SS$ be the set of unordered $mn$-tuples $\x$ of intersection
points between the lifts of $\alpha$ and $\beta$ circles such that
each such lift contains exactly one point of $\x$. (We will give a
more explicit characterization of the elements of $\SS$ later.) Let
$C$ be the $\Z/2$-vector space generated by $\SS$. Define a
differential $\partial$ on $C$ by making the coefficient of $\y$ in
$\tilde\partial\x$ nonzero if and only if the following conditions
hold:
\begin{itemize}
\item All but two of the points in $\x$ are also in $\y$.
\item There is a rectangle $R \in \tilde\RR$ whose lower-left and
upper-right corners are in $\x$, whose upper-left and lower-right
corners are in $\y$, and which does not contain any point of $\x$ in
its interior.
\end{itemize}
In Section \ref{sec:grid}, we shall define two gradings (Alexander
and Maslov) on $C$, as well as a decomposition of $C$ as a direct
sum of complexes corresponding to spin$^c$ structures on
$\Sigma^m(K)$. We shall prove:
\begin{thm}
The homology of the complex $(C, \partial)$ is isomorphic as a
bigraded group to $\HFK(\Sigma^m(K), \tilde K; \Z/2) \otimes
V^{\otimes n-1}$, where $V \cong \Z/2 \oplus \Z/2$ with generators
in bigradings $(0,0)$ and $(-1,-1)$.
\end{thm}

In Section 2, we review the construction of Heegaard Floer homology
for knots using multi-pointed Heegaard diagrams. In Section 3, we
show how to obtain a Heegaard diagram for $(\Sigma^m(K),\tilde K)$
given one for $(S^3,K)$, and we use apply that discussion to grid
diagrams in Section 4, proving Theorem 1. In Section 5, we give the
values of $\HFK(\Sigma^m(K),\tilde K)$ for over fifty knots with up
to eleven crossings. (Grigsby \cite{G2} has shown how to compute
these groups for two-bridge knots, so our tables only include knots
that are not two-bridge.) Finally, we make some observations and
conjectures about these results in Section 6.

\medskip \noindent \textbf{Acknowledgments.} I am grateful to Peter
Ozsv\'ath for suggesting this problem, providing lots of guidance,
and reading a draft of this paper, and to John Baldwin, Tom Peters,
Josh Greene, Matthew Hedden, and especially Eli Grigsby for many
extremely helpful conversations.

\section{Review of Heegaard Floer homology for knots}

Let us briefly recall the basic construction of Heegaard Floer
homology for knots \cite{OSzKnot}. For simplicity, we work with
coefficients modulo 2. A \emph{multi-pointed Heegaard diagram} $\D =
(\Sigma, \bm\alpha, \bm\beta, \w, \z)$ consists of an oriented
surface $\Sigma$; two sets of closed, embedded curves
$\bm\alpha=\{\alpha_1, \dots, \alpha_{g+n-1}\}$ and
$\bm\beta=\{\beta_1, \dots, \beta_{g+n-1}\}$ (where $g = g(\Sigma)$
and $n\ge 1$), each of which spans a $g$-dimensional subspace of
$H_1(\Sigma;\Z)$; and two sets of basepoints, $\w = \{w_1, \dots,
w_n\}$ and $\z = \{z_1, \dots, z_n\}$, such that each component of
$\Sigma - \bigcup \alpha_i$ and each component of $\Sigma - \bigcup
\beta_i$ contains exactly one point of $\w$ and one point of $\z$.
We obtain an oriented 3-manifold $Y$ and a handlebody decomposition
$Y = H_\alpha \cup_\Sigma H_\beta$ by attaching 2-handles to $\Sigma
\times I$ along the circles $\alpha_i \times \{0\}$ and $\beta_i
\times \{1\}$ and then canonically filling in 3-balls. To obtain a
knot or link $K$, we connect the $w$ (resp.~$z$) basepoints to the
$z$ (resp.~$w$) basepoints with arcs in the complement of the
$\alpha$ (resp.~$\beta$) curves and push those arcs into $H_\alpha$
(resp.~$H_\beta$). The orientations are such that $K$ intersects
$\Sigma$ positively at the $z$ basepoints (where it is passing from
$H_\alpha$ to $H_\beta$) and negatively at the $w$ basepoints (where
it is passing from $H_\beta$ to $H_\alpha$).

In terms of Morse theory, we obtain a Heegaard diagram for a given
pair $(Y,K)$ by taking a self-indexing Morse function $f$ on $Y$ and
a Riemannian metric such that $K$ is a union of gradient flowlines
connecting all the index-0 and index-3 basepoints. We then define
$\Sigma$ as $f^{-1}(\frac32)$, the $\alpha$ (resp.~$\beta$) circles
as the intersections of $\Sigma$ with the ascending (resp.
descending) manifolds of index-1 (resp.~index-2) critical points of
$f$, and the $w$ (resp.~$z$) basepoints as the intersections of
$\Sigma$ with the segments of $K$ that go from the index-3 (resp.
index-0) critical points to the index-0 (resp.~index-3) critical
points. We then have $H_\alpha = f^{-1}([0,\frac32])$ and $H_\beta =
f^{-1}([\frac32,3])$.

The Heegaard Floer complex $\CFK(\D)$ is defined as follows. Let
$\T_\alpha$ and $\T_\beta$ be the images of $\alpha_1 \times \dots
\times \alpha_{g+n-1}$ and $\beta_1 \times \dots \times
\beta_{g+n-1}$ in the symmetric product $\Sym^{g+n-1}(\Sigma)$;
these are both embedded copies of $T^{g+n-1}$. The group $\CFK(\D)$
is the $\Z/2$-vector space generated by the (finitely many)
intersection points in $\T_\alpha \cap \T_\beta$, and the
differential $\partial$ is defined by taking counts of holomorphic
disks connecting intersection points:
\[\partial \x = \sum_{\y \in \T_\alpha \cap \T_\beta}
\sum_{\substack{ \phi \in \pi_2(\x,\y) \mid  \\ \mu(\phi)=1 \\
n_w(\phi) = n_z(\phi) = 0 }} \# \left( \frac{ \mathcal{M}(\phi)
}{\R} \right) \y.\] Each homotopy class of Whitney disks $\phi \in
\pi_2(\x,\y)$ has an associated \emph{domain} in $\Sigma$: a 2-chain
$D = \sum a_i D_i$, where the $D_i$ are components of $\Sigma -
\bigcup \alpha_i - \bigcup \beta_i$ (\emph{elementary domains}),
such that $\partial D$ is made of arcs of $\alpha$ curves that
connect each point of $\x$ to a point of $\y$ and arcs of $\beta$
curves that connect each point of $\y$ to a point of $\x$. Then
$n_\w(\phi)$ and $n_\z(\phi)$ are the multiplicities of the
elementary domains containing points of $\w$ and $\z$, respectively.
The \emph{Maslov index} $\mu(\phi)$ can be computed using a formula
due to Lipshitz \cite{L}:
\[\mu(\phi) = \sum_i a_i e(D_i) + p_\x(D) + p_\y(D),\]
where $p_\x(D)$ (resp.~$p_\y(D)$) equals the sum of the average of
the multiplicities of the domains at the four corners of each point
of $\x$ (resp.~$\y$), and $e(D_i)$, the \emph{Euler measure}, equals
$1-\frac{k}{2}$ when $D_i$ is a convex $2k$-gon. The coefficient of
$\y$ represents the number of holomorphic representatives of $\phi$
and generally depends on the choice of almost complex structure on
$\Sigma$. For suitable choices, the homology of the complex is then
isomorphic to $\HFK(Y,K) \otimes V^{\otimes n-1}$, where $V \cong
\Z/2 \oplus \Z/2$ with generators in bigradings $(-1,-1)$ and
$(0,0)$, and $\HFK(Y,K)$ is an invariant of the knot type of $K
\subset Y$.

To define the spin$^c$ structure $\s_\w(\x)$ associated to a
generator $\x$, let $N_\x$ be the union of regular neighborhoods of
the closures of the gradient flowlines through the points of $\x$
and $\w$. (Flowlines through the former connect index-1 and index-2
critical points of $f$; those through the latter connect index-0 and
index-3 critical points.) The gradient vector field $\vec\nabla f$
is non-vanishing on $Y-N_\x$ and hence defines a spin$^c$ structure
(using Turaev's formulation of spin$^c$ structures as homology
classes of non-vanishing vector fields \cite{T}). Let $\CFK(\D, \s)
\subset \CFK(\D)$ be the subspace generated by the generators $\x$
with $\s_\w(\x)=\s$. To test whether two generators $\x$ and $\y$
are in the same spin$^c$ structure, let $\gamma_{\x,\y}$ be a
1-cycle obtained by connecting $\x$ to $\y$ along the $\alpha$
circles and $\y$ to $\x$ along the $\beta$ circles, and let
$\epsilon(\x,\y)$ be its image in
\[H_1(Y) \cong H_1(\Sigma)/\Span([\alpha_i], [\beta_i] \mid i=1,
\dots, g+n-1).\] Then $\x$ and $\y$ are in the same spin$^c$
structure if and only if $\epsilon(\x,\y)=0$. In particular, if $\y$
appears in the boundary of $\x$, then $\epsilon(\x,\y)=0$, so
$\CFK(\D,\s)$ is a subcomplex. The homology of each of these
summands does not depend on the choice of complex, so there is a
natural splitting
\[\HFK(Y,K) = \bigoplus_{\s \in \Spin^c(Y)} \HFK(Y,K,\s).\]

If $K$ is nulhomologous, the \emph{Alexander grading} on $\CFK(Y,K)$
is defined as follows. For each generator $\x$, let
$\srel_{\w,\z}(\x) \in \underline{\Spin^c}(Y,K) = \Spin^c(Y_0(K))$
be the spin$^c$ structure on the zero-surgery $Y_0(K)$ obtained by
extending $\s_\w(\x)|_{Y-N(K)}$ over $Y_0(K)$. Given a Seifert
surface $F$ for $K$, the Alexander grading of $\x$ is $A_F(\x) =
\frac12 \gen{ c_1(\srel_{\w,\z}(\x)), [\hat F]}$, where $\hat F$ is
the capped-off Seifert surface in $Y_0(K)$. The Alexander grading is
always independent of the choice of $F$ up to an additive constant
and completely independent when $Y$ is a rational homology sphere.
The relative Alexander grading between two generators $\x$ and $\y$,
$A(\x,\y) = A(\x) - A(\y)$, can also be given as the linking number
of $\gamma_{\x,\y}$ and $K$ (i.e., the intersection number of
$\gamma_{\x,\y}$ with $F$), or by the formula $A(\x,\y) = n_\z(D) -
n_\w(D)$ when $\x$ and $\y$ are in the same spin$^c$ structure and
$\D$ is any domain connecting $\x$ to $\y$. The latter formula shows
that the complex $\CFK(\D,\s)$ splits according to Alexander
gradings, and hence \[\HFK(Y,K,\s) = \bigoplus_{i \in \Z}
\HFK(Y,K,\s,i).\]

If $\s \in \Spin^c(Y)$ is a torsion spin$^c$ structure, the
\emph{relative Maslov grading} between two generators $\x$ and $\y$
with $\s_\w(\x) = \s_\w(\y)=\s$ is given by $M(\x,\y) = \mu(D) -
2n_\w(D)$, where $D$ is any domain connecting $\x$ to $\y$. An easy
way to compute the relative Maslov grading between two generators in
the same spin$^c$ structure is to find a linear combination of
$\alpha$ and $\beta$ circles that is homologous to $\gamma_{\x,\y}$
(which is possible since $\gamma_{\x,\y} \equiv 0$ in $H_1(Y)$).
Then $\gamma_{\x,\y}$ minus this linear combination bounds a domain
$D$ in $\Sigma$ connecting $\x$ to $\y$, and we then apply
Lipshitz's formula to compute $\mu(D)$.

Moreover, if $Y$ is a rational homology sphere, the relative
$\Z$-gradings on the $\CFK(Y,K,\s)$ lift to an absolute $\Q$-grading
on all of $\CFK(Y,K)$. Lipshitz and Lee \cite{LL} show that it is
easy to compute the relative $\Q$-grading between two generators
that are not necessarily in the same spin$^c$ structure. Since
$H_1(Y)$ is finite, there exists $m\ge 1$ such that
$m\gamma_{\x,\y}$ is homologous to a linear combination of $\alpha$
and $\beta$ circles, so $m\gamma_{\x,\y}$ minus this combination
bounds a domain $D$. The relative Maslov $\Q$-grading between $\x$
and $\y$ is then $M(\x,\y) = \frac1m(\mu(D) - 2n_\w(D))$. The
absolute $\Q$-grading is more complicated, and we shall not discuss
it in this paper.

Call a diagram $\D$ \emph{good} if every elementary domain that does
not contain a basepoint is either a bigon or a square. Manolescu,
Ozsv\'ath, and Sarkar showed that in any good diagram, the
coefficient of $\y$ in $\partial\x$ is nonzero in two cases:
\begin{itemize}
\item All but one of the points of $\y$ are also in $\x$, and
the remaining two points are the vertices of a bigon without a
basepoint or a point of $\x$ in its interior.
\item All but two of the points of $\y$ are also in $\x$, and the
remaining four points are the vertices of a rectangle without a
basepoint or a point of $\x$ in its interior.
\end{itemize}
It follows that when $\D$ is a good diagram, the boundary map can be
determined simply from the combinatorics of the diagram, without
reference to the choice of complex structure on $\Sigma$, so
$\HFK(Y,K)$ can be computed algorithmically.

If $K$ is a knot in $S^3$, then a grid diagram for $K$, drawn on a
torus as in Section 1, yields a Heegaard diagram $\D = (T^2,
\bm\alpha, \bm\beta, \w,\z)$ for the pair $(S^3,K)$, where the
$\alpha$ circles are the horizontal lines of the grid, the $\beta$
circles are the vertical lines, and the $w$ and $z$ basepoints are
the points marked $O$ and $X$, respectively. Every region of this
diagram is a square, so $\HFK(S^3,K)$ can be computed
combinatorially as above. Specifically, the generators correspond to
permutations of the set $\{1, \dots, n\}$, and the Alexander and
Maslov gradings of each generator can be given by simple formulae
(discussed later). Using this diagram, Baldwin and Gillam \cite{BG}
have computed $\HFK(S^3,K)$ for all knots with up to 12 crossings.
Additionally, Manolescu, Ozsv\'ath, Szab\'o, and Thurston
\cite{MOSzT} give a self-contained proof that this construction
yields a knot invariant. (See also Sarkar and Wang \cite{SW}, who
show how to obtain good diagrams for knots in arbitrary
3-manifolds.)

\section{Heegaard diagrams for cyclic branched covers of knots}

Given a knot $K \subset S^3$ and an integer $m \ge 2$, there is a
well-known construction of a 3-manifold $\Sigma^m(K)$ and an
$m$-fold branched covering map $\pi: \Sigma^m(K) \to S^3$ whose
downstairs branch locus is $K$ and whose upstairs branch locus is a
knot $\tilde K \subset \Sigma^m(K)$. The manifold $\Sigma^m(K)$ can
be constructed from $m$ copies of $S^3 - \operatorname{int} F$,
where $F$ is a Seifert surface for $K$, by connecting the negative
side of a bicollar of $F$ in the $k\Th$ copy to the positive side in
the $(k+1)\Th$ (indices modulo $m$). The inverse image of $K$ in
$\Sigma^m(K)$ is a knot $\tilde K$, which is nulhomologous because
it bounds a Seifert surface (any of the lifts of the original
Seifert surface $F$). For the details of this construction, see
Rolfsen \cite{Ro}.

The group of covering transformations of $\Sigma^m(K)\to S^3$ is
cyclic of order $m$, generated by a map $\tau_m: \Sigma^m(K) \to
\Sigma^m(K)$ that takes the $k\Th$ copy of $S^3\minus
\operatorname{int} F$ to the $(k+1)\Th$ (indices modulo $m$). If
$\gamma$ is a $1$-cycle in $S^3$, then by using transfer
homomorphisms, we see that for any lift $\tilde\gamma$, the equation
\begin{equation} \label{eqn:transfer}
\sum_{k=0}^{m-1} \tau_{m*}^k(\tilde \gamma) = 0
\end{equation}
holds in $H_1(\Sigma^m(K);\Z)$. In particular, when $m=2$, we have
$\tau_{2*}(\tilde \gamma) = -\tilde \gamma$.

When $m$ is a power of a prime $p$, the group $H_1(\Sigma^m(K);\Z)$
is then finite and contains no $p^r$-torsion for any $r$
\cite[p.~16]{Gor}. The order of $H_1(\Sigma^m(K))$ is equal to
$\prod_{j=0}^{m-1} \Delta_K(\omega^j)$, where $\Delta_K$ is the
Alexander polynomial of $K$, and $\omega$ is a primitive $m\Th$ root
of unity \cite[p.~149]{F}. In particular, note that the action of
the deck transformation group on $H_1(\Sigma^m(K);\Z)$ has no
nonzero fixed points: if $\tau_{m*}(\alpha)=\alpha$, then \[0 =
\alpha + \tau_{m*}(\alpha) + \dots + \tau_{m*}^{m-1}(\alpha) =
m\alpha,\] by Equation \ref{eqn:transfer}, so $\alpha=0$.

Let $\D = (S, \bm\alpha, \bm\beta, \w, \z)$ be a multi-pointed
Heegaard diagram for $K \subset S^3$ with genus $g$ and $n$
basepoint pairs.\footnote{In the discussion that follows, we denote
the Heegaard surface by $S$ rather than $\Sigma$ to avoid confusion
with the notation $\Sigma^m(K)$.} If $f: S^3 \to \R$ is a
self-indexing Morse function compatible with $\D$, then $\tilde f =
f \circ \pi: \Sigma^m(K) \to \R$ is a self-indexing Morse function
for the pair $(\Sigma^m(K),\tilde K)$ whose critical points are
simply the inverse images of the critical points of $f$. This
function induces a Heegaard splitting $\Sigma^m(K) = \tilde H_\alpha
\cup_{\tilde S} \tilde H_\beta$ that projects onto the Heegaard
splitting of $S^3$. A simple Euler characteristic argument shows
that the genus of the new Heegaard surface $\tilde S = \pi^{-1}(S)$
is $h=mg + (m-1)(n-1)$. Each $\alpha$ and $\beta$ circle in $S$
bounds a disk in $S^3 \minus K$ and hence has $m$ distinct preimages
in $\Sigma^m(K)$. Thus, we obtain a Heegaard diagram $\tilde\D =
(\tilde S, \tilde{\bm\alpha}, \tilde{\bm\beta}, \tilde\w,
\tilde\z)$, where $\tilde S$ is a surface of genus $h$ and
$\tilde{\bm\alpha}, \tilde{\bm\beta}, \tilde\w, \tilde\z$ are the
inverse images of the corresponding objects under the covering map.

We may arrange that the Heegaard surface $F$ intersects $S$ in $n$
arcs, each connecting a $z$ basepoint to a $w$ basepoint. Note that
each $\alpha$ or $\beta$ circle intersects $F$ algebraically zero
times, since, e.g., $\alpha_i \cdot F = \lk(\alpha_i, K) = K \cdot
D_{\alpha_i} = 0$, where $D_{\alpha_i}$ is a spanning disk for
$\alpha_i$. To obtain the diagram $\tilde\D$ directly, we may
connect $m$ copies of $\D$ by using the arcs of $F \cap S$ as branch
cuts. A complex structure on $S$ naturally yields a complex
structure on $\tilde S$ that makes the projection $\pi: \tilde S \to
S$ and the covering transformation $\tau_m: \tilde S \to \tilde S$
holomorphic.

The generators of the complex $\CFK(\tilde\D)$ may be described as
follows:

\begin{lemma} \label{decompose}
Any generator $\x$ of $\CFK(\tilde \D)$ can be decomposed
(non-uniquely) as $\x = \tilde\x_1 \cup \dots \cup \tilde \x_m$,
where $\x_1, \dots, \x_m$ are generators of $\CFK(\D)$, and
$\tilde\x_i$ is a lift of $\x_i$ to $\tilde \D$.
\end{lemma}

\begin{proof}
Given a generator $\x$ of $\CFK(\tilde\D)$, let $\bar\x$ be its
image under the natural map $\Sym^{mn}(\tilde S) \to \Sym^{mn}(S)$,
consisting of $mn$ points of $\Sigma$ (possibly repeated) such that
each $\alpha$ circle and each $\beta$ circle contains exactly $m$
points. It is then easy to partition $\bar\x$ into $m$ subsets
$\x_1, \dots, \x_m$, each of which is a generator of $\CFK(\D)$ as
required. Note that this choice of partition is not unique.
\end{proof}

Given a generator $\x_0$ of $\CFK(\D)$, let $L(\x_0)$ denote the
generator of $\CFK(\tilde\D)$ consisting of all $m$ lifts of each
point of $\x_0$. Using the action of the deck transformation
$\tau_m$ on $\D$, we may write $L(\x_0) = \tilde \x_0 \cup
\tau_m(\tilde\x_0) \cup \dots \cup \tau_m^{m-1}(\tilde\x)$, where
$\tilde\x_0$ is any lift of $\x_0$ to $\tilde D$.

\begin{lemma} \label{canonical}
All generators of $\CFK(\tilde\D)$ of the form $\x = L(\x_0)$ are in
the same spin$^c$ structure, denoted $\s_0$ and called the
\emph{canonical spin$^c$ structure} on $\Sigma^m(K)$.
\end{lemma}

\begin{proof}
(Adapted from Grigsby \cite{G1}.) Let $\x_0$ and $\y_0$ be
generators of $\CFK(\D)$; we shall show that $L(\x_0)$ and $L(\y_0)$
are in the same spin$^c$ structure. Let $\gamma_{\x_0,\y_0}$ be a
1-cycle joining $\x_0$ and $\y_0$ as above, and let $\tilde
\gamma_{\x_0,\y_0}$ be a lift of $\gamma_{\x_0,\y_0}$ to $\tilde S$.
Then the 1-cycle \[\tilde \gamma_{\x_0,\y_0} + \tau_{m*}(\tilde
\gamma_{\x_0,\y_0}) + \dots + \tau_{m*}^{m-1}( \tilde
\gamma_{\x_0,\y_0})\] connects $L(\x_0)$ and $L(\y_0)$. Then
$\epsilon(L(\x_0),L(\y_0)) =0$ by Equation \ref{eqn:transfer}, so
$L(\x_0)$ and $L(\y_0)$ are in the same spin$^c$ structure.
\end{proof}

\begin{remark}
When $K$ is a two-bridge knot and $m=2$, Grigsby shows that for a
specific diagram $\D$, the map $L$ extends to an isomorphism of
bigraded chain complexes $\CFK(\D) \to \CFK(\tilde D, \s_0)$.
Therefore, for any two-bridge knot $K$, $\HFK( \Sigma^2(K), \tilde
K, \s_0) \cong \HFK(S^3, K)$. In general, though, $L$ is not even a
chain map.
\end{remark}

The spin$^c$ structure $\s_0$ often also admits a more intrinsic
characterization. Assume $m$ is a prime power. If $f: S^3 \to \R$ is
a self-indexing Morse function for $(S^3,K)$ as above, then its
pullback $\tilde f: \Sigma^m(K) \to \R$ is $\tau_m$-invariant. Using
a Riemannian metric on $\Sigma^m(K)$ that is the pullback of a
metric on $S^3$, the gradient $\vec\nabla \tilde f$ is
$\tau_m$-invariant and projects onto $\vec\nabla f$, and the
flowlines for $\tilde f$ are precisely the lifts of flowlines for
$f$. If $N_{\x_0}$ is the union of neighborhoods of flowlines
through the points of $\x_0$ and $\w$, where $\x_0$ is a generator
of $\CFK(\D)$, then $\pi^{-1}(N_{\x_0})$ is the union of
neighborhoods of flowlines through the points of $L(\x_0)$ and can
be denoted $N_{L(\x_0)}$ as in Section 2. By suitably modifying
$\vec\nabla \tilde f$ on $N_{L(\x_0)}$, we may obtain a
$\tau_m$-invariant vector field that determines
$\s_{\tilde\w}(L(\x_0))=\s_0$. It follows that $\s_0$ is fixed under
the action of $\tau_m$ on $\Spin^c(\Sigma^m(K))$.\footnote{In
general, spin$^c$ structures can always be pulled back under a local
diffeomorphism using the vector field interpretation. Specifically,
if $F:M \to N$ is a local diffeomorphism and $\xi$ is a nonvanishing
vector field on $N$ that determines a given spin$^c$ structure $\s
\in \Spin^c(N)$, then $F^*(\s) \in \Spin^c(\Sigma^m(K)_0)$ is
determined by the vector field $(F_*)^{-1}(\xi)$. The first Chern
class is natural under this pullback.} Now, if $\s_0'$ is another
spin$^c$ structure fixed under the action of $\tau_m$, then the
difference between $\s_0$ and $\s_0'$ is a class in
$H_1(\Sigma^m(K);\Z)$ that is fixed by $\tau_m$, hence equals zero.
Thus, $\s_0$ is uniquely characterized by the property that
$\tau_m^*(\s_0) = \s_0$. For more about the significance of $\s_0$,
see \cite{GRS}.

We now consider the Alexander gradings in $\CFK(\tilde\D)$.

\begin{prop} \label{alexavg}
If $\x = \tilde\x_1 \cup \dots \cup \tilde\x_m$ as in Lemma
\ref{decompose}, then the Alexander grading of $\x$ (computed with
respect to a Seifert surface for $\tilde K$ that is a lift of a
Seifert surface for $K$) is equal to the average of the Alexander
gradings of $\x_1, \dots, \x_m$.
\end{prop}

\begin{proof}
We first consider the relative Alexander gradings. Let $F \subset
S^3$ be a Seifert surface for $K$, and let $\tilde F$ be a lift of
$F$ to $\Sigma^m(K)$. The translates $\tilde F, \tau_m(\tilde F),
\dots, \tau_m^{m-1}(\tilde F)$ are all Seifert surfaces for $\tilde
K$. The relative Alexander grading between two generators does not
depend on the choice of Seifert surface, so for generators $\x, \y$
of $\CFK(\tilde D)$, we have
\[m A (\x,\y) = \gamma_{\x,\y} \cdot \tilde F + \gamma_{\x,\y} \cdot
\tau_m(\tilde F) + \dots + \gamma_{\x,\y} \cdot \tau_m^{m-1} (\tilde
F),\] where $\gamma_{\x,\y}$ is a 1-cycle connecting $\x$ and $\y$
as above. The projection $\pi_*(\gamma_{\x,\y})$ is a 1-cycle in $S$
that goes from $\bar \x$ to $\bar \y$ along $\alpha$ circles and
from $\bar \y$ to $\bar \x$ along $\beta$ circles. Every
intersection point of $\gamma_{\x,\y}$ with one of the lifts of $F$
corresponds to an intersection point of $\pi_*(\gamma_{\x,\y})$ with
$F$, so
\[\gamma_{\x,\y} \cdot \tilde F + \gamma_{\x,\y} \cdot
\tau_m(\tilde F) + \dots + \gamma_{\x,\y} \cdot \tau_m^{m-1} (\tilde
F) = \pi_*(\gamma_{\x,\y}) \cdot F.\] The restriction of
$\pi_*(\gamma_{\x,\y})$ to any $\alpha$ or $\beta$ circle consists
of $m$ (possibly constant or overlapping) arcs. By perhaps adding
copies of the $\alpha$ or $\beta$ circle, we can arrange that these
arcs connect a point of $\x_1$ with a point of $\y_1$, a point of
$\x_2$ with a point of $\y_2$, and so on. In other words,
\[\pi_*(\gamma_{\x,\y}) \equiv \gamma_{\x_1,\y_1} + \dots +
\gamma_{\x_m,\y_m}\] modulo the $\alpha$ and $\beta$ circles in
$\D$, whose intersection numbers with $F$ are zero. Therefore,
\[
\begin{split}
A(\x,\y) &= \frac1m (\gamma_{\x_1,\y_1} + \dots, +
\gamma_{\x_m,\y_m}) \cdot F \\
&= \frac1m (A(\x_1,\y_1) + \dots + A(\x_m,\y_m)).
\end{split}
\]
Thus, the Alexander grading of a generator of $\CFK(\tilde\D)$ is
given up to an additive constant by the average Alexander grading of
its parts.

To pin down the additive constant, first note that the branched
covering map $\pi: \Sigma^m(K) \to S^3$ extends to an
\emph{unbranched} covering map from the zero-surgery on $\tilde K$
to the zero-surgery on $K$, $\pi_0: \Sigma^m(K)_0 \to S^3_0$. Since
this is a local diffeomorphism, it is possible to pull back spin$^c$
structures. Let $\x_0$ be a generator $\CFK(\D)$ in Alexander
grading 0, and let $\x = L(\x_0)$. As in the discussion following
Lemma \ref{canonical}, we may find a nonvanishing vector field that
determines $\s_{\tilde\w}(\x) = \s_0$ and is $\tau_m$-equivariant.
The unique extension (up to isotopy) of this vector field to
$\Sigma^m(K)_0$ can also be made $\tau_m$-invariant, so it is the
pullback of an extension to $S^3_0$ of a vector field determining
$\s_\w(\x_0)$. It follows that $\srel_{\tilde\w,\tilde\z}(\x) =
\pi_0^* (\srel_{\w,\z}(\x_0))$. Now, if $\hat{\tilde F} \subset
Y_0(\tilde K)$ is obtained by capping off $\tilde F$ in the
zero-surgery, then $\pi_{0*} [\hat{\tilde F}] = [\hat F]$ in
$H_2(S_3^0;\Z)$. Therefore,
\begin{equation*}
\begin{split}
A_{\tilde F}(\x) &= \frac12 \gen{ c_1(\srel_{\tilde\w,
\tilde\z}(\x)),
[ \hat{\tilde F} ] } \\
&= \frac12 \gen { c_1(\pi_0^*(\srel_{\w,\z}(\x_0))), [ \hat{\tilde F} ] } \\
&= \frac12 \gen { c_1(\srel_{\w,\z}(\x_0)), \pi_{0*} [ \hat{\tilde F} ] } \\
&= \frac12 \gen { c_1(\srel_{\w,\z}(\x_0)), [ \hat F ] } \\
&= 0 = A_F(\x_0).
\end{split}
\end{equation*}
Thus, the additive constant $C$ must equal 0.
\end{proof}

Next, we consider the domains in $\tilde D$. Any simply-connected
elementary domain $D$ of $\D$ that does not contain a basepoint is
evenly covered, so its preimage in $\tilde\D$ consists of $m$
disjoint domains each diffeomorphic to $D$. On the other hand, a
domain containing exactly one basepoint is covered by a single
connected domain with $m$ times as many sides as the original one.
In particular, if $\D$ is a good diagram, then $\tilde\D$ is also
good. It follows that the domains that count for the boundary in
$\CFK(\tilde\D)$ are precisely the lifts of the domains that count
for the boundary of $\D$.

We conclude with a few comments about the symmetries in the case
where $m=2$. The order of $H_1(\Sigma^2(K);\Z)$ is equal to the
determinant of $K$, $\det K = \Delta_K(-1)$, which is always odd. As
mentioned above, the non-trivial deck transformation $\tau_2$ acts
on $H_1(\Sigma^2(K);\Z)$ by multiplication by $-1$. The set
$\Spin^c(\Sigma^2(K))$ of spin$^c$ structures on $\Sigma^2(K)$ is an
affine set for $H_1(\Sigma^2(K);\Z)$ and can be identified with the
latter by sending the canonical spin$^c$ structure $\s_0$ to zero.
Under this identification, both conjugation ($\s \mapsto \bar\s$)
and pullback under $\tau_2$ ($\s \mapsto \tau_2^*(\s)$) are given by
with multiplication by $-1$, so $\tau_2^*(\s)=\bar\s$. Since the
diagram $\tilde\D$ is $\tau_2$-equivariant, $\tau_2$ induces an
isomorphism of bigraded groups \[\HFK(\Sigma^2(K), \tilde K, \s) \to
\HFK(\Sigma^2(K), \tilde K, \bar\s).\] On the other hand, it is a
standard fact \cite[Prop. 3.10]{OSzKnot} that \[\HFK_j(Y,K,\s,i)
\cong \HFK_{j-2i} (Y,K,\bar\s,-i).\] Therefore, to compute
$\HFK(\Sigma^2(K), \tilde K)$, it suffices to consider only one out
of every pair of conjugate, non-canonical spin$^c$ structures, and
to consider only the generators that lie in non-negative Alexander
grading. Additionally, note that since $\Sigma^2(K)$ is a rational
homology sphere, the Maslov $\Z$-grading lifts to a $\Q$-grading
that extends across all spin$^c$ structures.

\section{Grid diagrams and cyclic branched covers} \label{sec:grid}

As described in Section 1, any oriented knot $K \subset S^3$ can be
represented by means of a grid diagram. By drawing the grid diagram
on a standardly embedded torus in $S^3$, we may think of the grid
diagram as a genus 1, multi-pointed Heegaard diagram $\D=(T^2,
\bm\alpha, \bm\beta, \w, \z)$ for the pair $(S^3,K)$, where the
$\alpha$ circles are the horizontal lines of the grid, the $\beta$
circles are the vertical lines, the $w$ basepoints are in the
regions marked $O$, and the $z$ basepoints are in the regions marked
$X$.

We label the $\alpha$ circles $\alpha_0, \dots, \alpha_{n-1}$ from
bottom to top and the $\beta$ circles $\beta_0, \dots, \beta_{n-1}$
from left to right. Each $\alpha$ circle intersects each $\beta$
circle exactly once: $\beta_i \cap \alpha_j = \{x_{ij}\}$.
Generators of the Heegaard Floer chain complex $\CFK(\D)$ then
correspond to permutations of the index set $\{0, \dots, n-1\}$ via
the correspondence $\sigma \mapsto (x_{0,\sigma(0)}, \dots,
x_{n-1,\sigma(n-1)})$. The diagram is good, so the differential can
be computed combinatorially as described in Section 2. Specifically,
the coefficient of $\y$ in $\partial\x$ is 1 if all but two of the
points of $\x$ and $\y$ agree and there is a rectangle embedded in
the torus with points of $\x$ as its lower-left and upper-right
corners, points of $\y$ as its lower-right and upper-left corners,
and no basepoints or points of $\x$ in its interior, and 0
otherwise.

For each grid point $x$, let $w(x)$ denote the winding number of the
knot projection around $x$. Let $p_1, \dots, p_{8n}$ (repetitions
allowed) denote the vertices of the $2n$ squares containing
basepoints, and set
\[a = \frac{1-n}{2} + \frac18 \sum_{i=1}^{8n} w(p_i).\] According to
Manolescu, Ozsv\'ath, and Sarkar \cite{MOS}, the Alexander grading
of a generator $\x$ of $\CFK(\D)$ is given by the formula
\begin{equation} \label{eqn:alexdownstairs}
A(\x) = a - \sum_{x\in\x} w(x).
\end{equation}
There is also a formula for the Maslov grading of a generator, but
it is not relevant for our purposes.

A Seifert surface for $K$ may be seen as follows. Isotope $K$ so
that it lies entirely within $H_\alpha$ by letting the arcs of $K
\cap H_\beta$ fall onto the boundary torus. In fact, it lies within
a ball contained in $H_\alpha$ since the knot projection in the grid
diagram never passes through the left edge of the grid. Take a
Seifert surface $F$ contained in this ball, and then isotope $F$ and
$K$ so that $K$ returns to its original position. $F$ then
intersects the Heegaard surface $T^2$ in $n$ arcs, one connecting
the two basepoints in each column of the grid diagram, and it
intersects $H_\beta$ in strips that lie above these arcs. The
orientations of $K$ and $S^3$ imply that the positive side of a
bicollar for $F$ lies on the \emph{right} of one of these strips
when the $X$ is above the $O$ and on the \emph{left} when the $O$ is
above the $X$.

By the results of Section 3, it follows that $\tilde\D = (\tilde T,
\tilde{\bm\alpha}, \tilde{\bm\beta}, \tilde w, \tilde z)$, where
$\tilde T$ is the surface defined in Section 1 and
$\tilde{\bm\alpha}, \tilde{\bm\beta}, \tilde w, \tilde z$ are the
lifts of the corresponding objects in $\D$, is a good Heegaard
diagram for $(\Sigma^m(K), \tilde K)$.

For computational purposes, the generators of $\CFK(\tilde\D)$ can
be described easily as follows. For any $i=0, \dots, n-1$ and $j=0,
\dots, n-1$, each lift of $\beta_i$ meets exactly one lift of
$\alpha_j$. Specifically, let $\tilde\beta_j^k$ denote the lift of
$\beta_j$ on the $k\Th$ copy of $\D$ (for $k=0, \dots, m-1$). Let
$\tilde\alpha_j^k$ denote the lift of $\alpha_j$ that intersects the
leftmost edge of the $k\Th$ grid diagram ($\tilde\beta_0^k$). Let
$\tilde x_{i,j}^k$ denote the lift of $x_{i,j}$ on the $k\Th$
diagram. Define a map $g: \Z/n \times \Z/n \times \Z/m \to \Z/m$ by
$g(i,j,k) = k - w(x_{i,j}) \bmod m$. The lift of $\alpha_j$ that
meets a particular $\tilde\beta_i^k$ is given by the following
lemma:

\begin{lemma}
The point $\tilde x_{i,j}^k$ is the intersection between
$\tilde\beta_i^k$ and $\tilde\alpha_j^{g(i,j,k)}$.
\end{lemma}

\begin{proof}
We induct on $i$. For $i=0$, we have $w(x_{0,j})=0$, and by
construction $\tilde\alpha_j^k$ meets $\tilde\beta_0^k$. For the
induction step, let $\overrightarrow{x_{i,j} x_{i+1,j}}$ be the
segment of $\alpha_j$ from $x_{i,j}$ to $x_{i+1,j}$. Note that
$w(x_{i+1,j})$ is equal to $w(x_{i,j}) + 1$ if
$\overrightarrow{x_{i,j} x_{i+1,j}}$ passes below the $X$ and above
the $O$ in its column, $w(x_{i,j})-1$ if it passes above $X$ and
below $O$, and $w(x_{i,j})$ otherwise. Similarly, if $\tilde
x_{i,j}^k$ lies on $\tilde\alpha_j^l$, then by the previous
discussion, $\tilde x_{i+1,j}^k$ lies on $\tilde\alpha_j^{l-1}$ in
the first case, on $\tilde\alpha_j^{l+1}$ in the second, and on
$\tilde\alpha_j^l$ in the third (upper indices modulo $m$). This
proves the induction step.
\end{proof}

We may then identify the generators of $\CFK(\tilde\D)$ with the set
of $m$-to-one maps
\[\phi: \{0, \dots, n-1\} \times \{0, \dots, m-1\} \to \{0, \dots,
n-1\}\] such that for each $j=0, \dots, n-1$, the function $g(\cdot,
j, \cdot)$ assumes all $m$ possible values on $\phi^{-1}(j)$. In
other words, if we shade the $m$ lifts of each $\alpha$ with
different colors as in Figure \ref{fig:differential} and arrange the
copies of $T$ horizontally, a generator is a selection of $mn$ grid
points so each column contains one point and each row contains $m$
points, one of each color. It is not difficult to enumerate such
maps algorithmically.

The differentials in $\CFK(\tilde D)$ are easy to compute. Since all
of the regions of $\tilde D$ that do not contain basepoints are
rectangles, the only domains that count for the differential are
rectangles, as described above. These are precisely the lifts of the
domains in $\D$ that count for the differential of $\CFK(\D)$. This
proves Theorem 1.

To compute the Alexander grading of a generator $\x$, we decompose
it into $\tilde \x_1 \cup \dots \cup \tilde \x_m$ using Lemma
\ref{decompose}, and then use Proposition \ref{alexavg} and Equation
\ref{eqn:alexdownstairs} to write:
\begin{equation}
\begin{split}
A_{\tilde F} (\x) &= \frac1m ( A_{F}(\x_1) +
\dots + A_{F}(\x_m) ) \nonumber \\
&= \frac1m \sum_{k=1}^m \left( a-\sum_{x \in \x_k} w(x) \right) \nonumber \\
&= a - \frac1m \sum_{k=1}^m \sum_{x \in \tilde\x_k} w(\pi(x)) \nonumber \\
&= a - \frac1m \sum_{x \in \x} w(\pi(x)). \label{alexupstairs}
\end{split}
\end{equation}

To split up the generators of $\CFK(\tilde\D)$ according to spin$^c$
structures, we simply need to be able to express $\epsilon(\x,\y)$
in terms of a presentation $H_1(\Sigma^m(K);\Z)$. Since
\[H_1(\Sigma^m(K);\Z) \cong H_1(\tilde T) / \Span
([\tilde\alpha_i^k], [\tilde\beta_i^k] \mid i \in \Z/n, k \in \Z/m
),\] we can obtain such a presentation by taking a basis for
$H_1(\tilde S)$ and imposing relations obtained by expressing
$\tilde\alpha$ and $\tilde\beta$ curves in terms of that basis.

\begin{figure} \label{fig:h1sigma}
\includegraphics[scale=0.75]{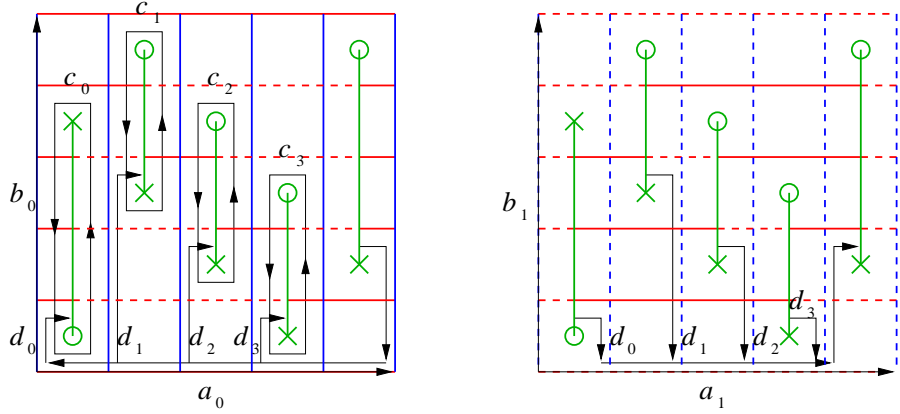}
\centering \caption{A symplectic basis for $H_1(\tilde T;\Z)$.}
\end{figure}

In the case where $m=2$, we may view $\tilde T$ as the union of two
$n$-times-punctured tori $T_0, T_1$, glued along their boundaries.
It is then easy to write down a symplectic basis for $H_1(\tilde
T;\Z)$. Specifically, let $(a_i, b_i)$ ($i=0,1$) be the standard
basis for $H_1(T_i;\Z)$, where $a_i$ is the bottom edge of the grid
diagram (oriented to the right) and $b_i$ is the left edge (oriented
upwards), so that $a_i \cdot b_i = 1$. Let $c_j$ ($j=0,\dots, n-2$)
be a loop in $T_1$ that goes once counterclockwise around the $j\Th$
branch cut (counted from the left), and let $d_j$ be a loop that
passes from the right side of the $(n-1)\Th$ branch cut to the left
side of the $j\Th$ branch cut in $T_0$ and from the right side of
the $j\Th$ branch cut to the left side of the $(n-1)\Th$ branch cut
in $T_1$, passing below all of the other branch cuts. (See Figure
\ref{fig:h1sigma}.) Then $c_j \cdot d_j=1$, and all other
intersection numbers are zero. It is not hard to see that the $a_i$,
$b_i$, and $c_j$ are all killed in $H_1(Y)$, and the remaining
relators are alternating sums of $d_j$ given by the
$\tilde\alpha_i^0$ circles. This presentation can then be reduced to
Smith normal form for easy use. For instance, in the right-handed
trefoil example shown in Figure \ref{fig:h1sigma},
\[
\begin{split}
H_1(\Sigma^2(K);\Z) &\cong \Z^4 \gen{d_0,\dots,d_3} / (d_0 - d_3,
d_0-d_2+d_3, d_0-d_1+d_2, d_1) \\
&\cong \Z/3.
\end{split}
\]
Computing $\epsilon(\x,\y)$ is then just a matter of counting how
many times a 1-cycle representative $\gamma_{\x,\y}$ passes through
the branch cuts, weighting the cuts appropriately.

The relative Maslov grading between two generators (an integer if
they are in the same spin$^c$ structure, and a rational number
otherwise) can be computed as described in Section 2. Because all
the basepoints in the Heegaard diagrams used in this paper are
contained in octagonal regions, it is not possible to compute the
absolute Maslov gradings or the spectral sequence from
$\HFK(\Sigma^2(K),\tilde K)$ to $\HF(\Sigma^2(K))$ combinatorially.
However, in many instances, the groups $\HF(\Sigma^2(K))$, or at
least the correction terms $d(\Sigma^2(K),\s)$, can be computed via
other means \cite{JN, OSzUnknot}. In such cases, it is often
possible to pin down the absolute Maslov gradings for
$\HFK(\Sigma^2(K), \tilde K)$. Specifically, the relative Maslov
$\Q$-grading and the action of $H_1(\Sigma^2(K))$ on
$\Spin^c(\Sigma^2(K))$ usually provide enough information to match
the groups $\HFK(\Sigma^2(K), \tilde K, \s)$ up with the rational
numbers $d(\Sigma^2(K),\s)$ that are computed via some other means.
If there is a spin$^c$ structure $\s$ in which $\HFK(\Sigma^2(K),
\tilde K,\s)$ has rank 1, then the absolute Maslov grading of that
group equals the corresponding $d$ invariant, and the rest of the
absolute gradings are completely determined.

\section{Results}

The tables that follow list the ranks for $\HFK(\Sigma^2(K),\tilde
K;\Z/2)$ by means of the Poincar\'e polynomials:
\[p_\s(q,t) = \sum_{i,j} \dim_{\Z/2} \HFK_j(\Sigma^2(K), \tilde K,
\s, i; \Z/2) t^i q^j.\] The Maslov $\Q$-gradings are normalized so
that in the canonical spin$^c$ structure $\s_0$, the nonzero
elements in Alexander grading $g(K)$ have Maslov grading $g(K)$. For
each knot, the first line gives $p_{\s_0}(q,t)$, and each subsequent
line gives $p_\s(q,t)$ for a pair of conjugate spin$^c$ structures.
We identify spin$^c$ structures with elements of
$H_1(\Sigma^2(K);\Z)$, which is either a cyclic group or the sum of
two cyclic groups, taking $\s_0$ to $0$. (Of course, the choice of
basis for $H_1(\Sigma^2(K);\Z)$ is not canonical.) In each spin$^c$
structure, most of the nonzero groups lie along a single diagonal;
the terms corresponding to the groups not on that diagonal are
underlined.

These results were computed using a program written in C++ and
\emph{Mathematica}, based on Baldwin and Gillam's program \cite{BG}
for computing $\HFK(S^3,K)$. Most of the grid diagrams were obtained
using Marc Culler's program \emph{Gridlink} \cite{MC}. Using
available computer resources, it was possible to compute
$\HFK(\Sigma^2(K),\tilde K)$ for all the three-bridge knots with up
to eleven crossings and arc index $\le 9$, and for many knots with
arc index 10. (Grigsby \cite{G2} has a much more efficient algorithm
for computing $\HFK(\Sigma^2(K),\tilde K)$ when $K$ is two-bridge,
so we do not list those knots here.)

\newpage
\begin{tiny}

\[
\begin{array}{cccl}
K & H_1(\Sigma^2(K);\Z) & \s & \sum_{i,j} \dim_{\Z/2}
\HFK_j(\Sigma^2(K), \tilde K, \s, i; \Z/2) t^i q^j \\
\\

8_5 &\ \Z/21 &0& q^{-3}t^{-3} + 3q^{-2}t^{-2} + 4q^{-1}t^{-1} + 5 +
4qt + 3q^2t^2 + q^3t^3 \\
&&\pm1& q^{5/21} (q^{-2}t^{-2} + 3q^{-1}t^{-1} + 3 + 3qt + q^2t^2) \\
&&\pm2& q^{20/21} \\
&&\pm3& q^{8/7} \\
&&\pm4& q^{17/21} (q^{-1}t^{-1} + 1 + qt) \\
&&\pm5& q^{20/21} \\
&&\pm6& q^{4/7} \\
&&\pm7& q^{2/3} (q^{-1}t^{-1} + 3 + qt) \\
&&\pm8& q^{5/21} (q^{-2}t^{-2} + 3q^{-1}t^{-1} + 3 + 3qt + q^2t^2) \\
&&\pm9& q^{2/7} (q^{-2}t^{-2} + 2q^{-1}t^{-1} + 3 + 2qt + q^2t^2) \\
&&\pm10& q^{17/21} (q^{-1}t^{-1} + 1 + qt) \\
\\

8_{10} &\Z/27 &0& q^{-3}t^{-3} + 3q^{-2}t^{-2} + 6q^{-1}t^{-1} + 7 +
6qt + 3q^2t^2 + q^3t^3 \\
&&\pm1& q^{7/27} (q^{-2}t^{-2} + 3q^{-1}t^{-1} + 5 + 3qt + q^2t^2) \\
&&\pm2& q^{1/27} \\
&&\pm3& q^{1/3} \\
&&\pm4& q^{4/27} (q^{-1}t^{-1} + 1 + qt) \\
&&\pm5& q^{13/27} \\
&&\pm6& q^{1/3} \\
&&\pm7& q^{-8/27} (q^{-1}t^{-1} + 1 + qt) \\
&&\pm8& q^{-11/27} (q^{-2}t^{-2} + q^{-1}t^{-1} + 1 + qt + q^2t^2) \\
&&\pm9& q^{-1}t^{-1} + 1 + qt \\
&&\pm10& q^{25/27} \\
&&\pm11& q^{10/27} (2q^{-1}t^{-1} + 3 + 2qt) \\
&&\pm12& q^{1/3} (q^{-2}t^{-2} + 2q^{-1}t^{-1} + 3 + 2qt + q^2t^2) \\
&&\pm13& q^{22/27} (q^{-1}t^{-1} + 1 + qt) \\
\\

8_{15} &\Z/33 &0& 3q^{-2}t^{-2} + 8q^{-1}t^{-1} + 11 + 8qt + 3q^2t^2
\\
&&\pm1& q^{13/33} (2q^{-1}t^{-1} + 3 + 2qt) \\
&&\pm2& q^{-14/33} (q^{-2}t^{-2} + q^{-1}t^{-1} + 1 + qt + q^2t^2) \\
&&\pm3& q^{6/11} \\
&&\pm4& q^{10/33} \\
&&\pm5& q^{-5/33} (q^{-1}t^{-1} + 1 + qt) \\
&&\pm6& q^{2/11} \\
&&\pm7& q^{10/33} \\
&&\pm8& q^{7/33} (q^{-1}t^{-1} + 1 + qt) \\
&&\pm9& q^{10/11} \\
&&\pm10& q^{13/33} (2q^{-1}t^{-1} + 3 + 2qt) \\
&&\pm11& q^{2/3} \\
&&\pm12& q^{-3/11} (q^{-1}t^{-1} + 1 + qt) \\
&&\pm13& q^{-14/33} (q^{-2}t^{-2} + q^{-1}t^{-1} + 1 + qt + q^2t^2) \\
&&\pm14& q^{7/33} (q^{-1}t^{-1} + 1 + qt) \\
&&\pm15& q^{-4/11} (q^{-2}t^{-2} + 2q^{-1}t^{-1} + 3 + 2qt + q^2t^2) \\
&&\pm16& q^{-5/33} (q^{-1}t^{-1} + 1 + qt) \\
\\

8_{16} &\Z/35 &0& q^{-3}t^{-3} + 4q^{-2}t^{-2} + 8q^{-1}t^{-1} + 9 +
8qt + 4q^2t^2 + q^3t^3 \\
&&\pm1& q^{16/35} (q^{-1}t^{-1} + 1 + qt) \\
&&\pm2& q^{29/35} \\
&&\pm3& q^{4/35} (q^{-1}t^{-1} + 1 + qt) \\
&&\pm4& q^{11/35} (q^{-2}t^{-2} + 3q^{-1}t^{-1} + 5 + 3qt + q^2t^2) \\
&&\pm5& q^{3/7} (q^{-1}t^{-1} + 3 + qt) \\
&&\pm6& q^{16/35} (q^{-1}t^{-1} + 1 + qt) \\
&&\pm7& q^{2/5} (2q^{-1}t^{-1} + 3 + 2qt) \\
&&\pm8& q^{9/35} \\
&&\pm9& q^{1/35} \\
&&\pm10& q^{5/7} (q^{-1}t^{-1} + 3 + qt) \\
&&\pm11& q^{11/35} (q^{-2}t^{-2} + 3q^{-1}t^{-1} + 5 + 3qt + q^2t^2) \\
&&\pm12& q^{29/35} \\
&&\pm13& q^{9/35} \\
&&\pm14& q^{3/5} \\
&&\pm15& q^{6/7} (q^{-1}t^{-1} + 1 + qt) \\
&&\pm16& q^{1/35} \\
&&\pm17& q^{4/35} (q^{-1}t^{-1} + 1 + qt) \\
\end{array}
\]

\[
\begin{array}{cccl}
K & H_1(\Sigma^2(K);\Z) & \s & \sum_{i,j} \dim_{\Z/2}
\HFK_j(\Sigma^2(K), \tilde K, \s, i; \Z/2) t^i q^j \\
\\

8_{17} &\Z/37 &0& q^{-3}t^{-3} + 4q^{-2}t^{-2} + 8q^{-1}t^{-1} + 11
+ 8qt + 4q^2t^2 + q^3t^3 \\
&&\pm1& q^{2/37} \\
&&\pm2& q^{8/37} (q^{-1}t^{-1} + 3 + qt) \\
&&\pm3& q^{18/37} \\
&&\pm4& q^{-5/37} (q^{-1}t^{-1} + 1 + qt) \\
&&\pm5& q^{13/37} (q^{-2}t^{-2} + 3q^{-1}t^{-1} + 3 + 3qt + q^2t^2) \\
&&\pm6& q^{-2/37} \\
&&\pm7& q^{-13/37} (q^{-2}t^{-2} + 3q^{-1}t^{-1} + 3 + 3qt + q^2t^2) \\
&&\pm8& q^{17/37} (q^{-1}t^{-1} + 1 + qt) \\
&&\pm9& q^{14/37} \\
&&\pm10& q^{-22/37} \\
&&\pm11& q^{-17/37} (q^{-1}t^{-1} + 1 + qt) \\
&&\pm12& q^{-8/37} (q^{-1}t^{-1} + 3 + qt) \\
&&\pm13& q^{5/37} (q^{-1}t^{-1} + 1 + qt) \\
&&\pm14& q^{22/37} \\
&&\pm15& q^{6/37} \\
&&\pm16& q^{-6/37} \\
&&\pm17& q^{-14/37} \\
&&\pm18& q^{-18/37} \\
\\

8_{18} &\Z/3 \oplus \Z/15 &(0,0)& q^{-3}t^{-3} + 5q^{-2}t^{-2} +
10q^{-1}t^{-1} + 13+ 10qt + 5q^2t^2 + q^3t^3 \\
&&\pm(0,1)& q^{7/15} (q^{-1}t^{-1} + 1 + qt) \\
&&\pm(0,2)& q^{-2/15} \\
&&\pm(0,3)& q^{1/5} (q^{-1}t^{-1} + 1 + qt) \\
&&\pm(0,4)& q^{7/15} (q^{-1}t^{-1} + 1 + qt) \\
&&\pm(0,5)& q^{-2/3} \\
&&\pm(0,6)& q^{-1/5} (q^{-1}t^{-1} + 1 + qt) \\
&&\pm(0,7)& q^{-2/15} \\
&&\pm(1,0)& q^{-2/3} \\
&&\pm(1,1)& q^{7/15} (q^{-1}t^{-1} + 1 + qt) \\
&&\pm(1,2)& q^{-7/15} (q^{-1}t^{-1} + 1 + qt) \\
&&\pm(1,3)& q^{-7/15} (q^{-1}t^{-1} + 1 + qt) \\
&&\pm(1,4)& q^{7/15} (q^{-1}t^{-1} + 1 + qt) \\
&&\pm(1,5)& q^{-2/3} \\
&&\pm(1,6)& q^{2/15} \\
&&\pm(1,7)& q^{-2/15} \\
&&\pm(1,8)& q^{-7/15} (q^{-1}t^{-1} + 1 + qt) \\
&&\pm(1,9)& q^{2/15} \\
&&\pm(1,10)& q^{2/3} \\
&&\pm(1,11)& q^{2/15} \\
&&\pm(1,12)& q^{-7/15} (q^{-1}t^{-1} + 1 + qt) \\
&&\pm(1,13)& q^{-2/15} \\
&&\pm(1,14)& q^{2/15} \\
\\

8_{19} &\Z/3 &0& q^{-3}t^{-3} + q^{-2}t^{-2} + \underline{q} + q^2t^2 + q^3t^3 \\
&&\pm1& q^{2/3}(q^{-1}t^{-1} + 1 + qt) \\
\\

8_{20} &\Z/9 &0& q^{-2}t^{-2} + 2q^{-1}t^{-1} + 3 + 2qt + q^2t^2 \\
&&\pm1& q^{7/9}(q^{-1}t^{-1} + 1 + qt) \\
&&\pm2& q^{1/9}(q^{-1}t^{-1} + 1 + qt) \\
&&\pm3& 1 \\
&&\pm4& q^{4/9} \\
\\

8_{21} &\Z/15 &1& q^{-2}t^{-2} + 4q^{-1}t^{-1} + 5 + 4qt + q^2t^2 \\
&&\pm1& q^{-2/15}(q^{-1}t^{-1} + 1 + qt) \\
&&\pm2& q^{7/15} \\
&&\pm3& q^{-1/5} (q^{-1}t^{-1} + 3 + qt) \\
&&\pm4& q^{-2/15}(q^{-1}t^{-1} + 1 + qt) \\
&&\pm5& q^{-1/3} \\
&&\pm6& q^{1/5} \\
&&\pm7& q^{7/15} \\
\end{array}
\]

\[
\begin{array}{cccl}
K & H_1(\Sigma^2(K);\Z) & \s & \sum_{i,j} \dim_{\Z/2}
\HFK_j(\Sigma^2(K), \tilde K, \s, i; \Z/2) t^i q^j \\
\\

9_{42} &\Z/7 &0& q^{-2}t^{-2} + 2q^{-1}t^{-1} + 2 + \underline{q} + 2qt + q^2t^2 \\
&&\pm1& q^{3/7} \\
&&\pm2& q^{5/7} (q^{-1}t^{-1} + 3 + qt) \\
&&\pm3& q^{6/7} (q^{-1}t^{-1} + 1 + qt) \\
\\

9_{43} &\Z/13 &0& q^{-3}t^{-3} + 3 q^{-2}t^{-2} + 2q^{-1}t^{-1} + 1
+ 2qt + 3q^2t^2  + q^3t^3 \\
&&\pm1& q^{10/13} (q^{-1}t^{-1} + 3 +qt) \\
&&\pm2& q^{1/13} (q^{-1}t^{-1} + 1 +qt) \\
&&\pm3& q^{12/13} \\
&&\pm4& q^{4/13} (q^{-2}t^{-2} + q^{-1}t^{-1} + 1 + qt + q^2t^2) \\
&&\pm5& q^{16/13} \\
&&\pm6& q^{9/13} (q^{-1}t^{-1} + 1 +qt) \\
\\

9_{44} &\Z/17 &0& q^{-2}t^{-2} + 4q^{-1}t^{-1} + 7 + 4qt + q^2t^2 \\
&&\pm1& q^{-8/17} \\
&&\pm2& q^{-15/17} (q^{-1}t^{-1} + 1 + qt) \\
&&\pm3& q^{-4/17} \\
&&\pm4& q^{8/17} \\
&&\pm5& q^{4/17} \\
&&\pm6& q^{-16/17} \\
&&\pm7& q^{-1/17} (q^{-1}t^{-1} + 1 + qt) \\
&&\pm8& q^{-2/17} (q^{-1}t^{-1} + 3 + qt) \\
\\

9_{45} &\Z/23 &0& q^{-2}t^{-2} + 6q^{-1}t^{-1} + 9 + 6qt + q^2t^2 \\
&&\pm1& q^{-8/23} (2q^{-1}t^{-1} + 3 + 2qt) \\
&&\pm2& q^{-9/23} \\
&&\pm3& q^{-3/23} (q^{-1}t^{-1} + 3 + qt) \\
&&\pm4& q^{-13/23} \\
&&\pm5& q^{7/23} \\
&&\pm6& q^{11/23} \\
&&\pm7& q^{-1/23} \\
&&\pm8& q^{-6/23} (q^{-1}t^{-1} + 1 + qt) \\
&&\pm9& q^{-4/23} (2q^{-1}t^{-1} + 3 + 2qt) \\
&&\pm10& q^{-18/23} (q^{-1}t^{-1} + 1 + qt) \\
&&\pm11& q^{-2/23} (q^{-1}t^{-1} + 1 + qt) \\
\\

9_{46} &\Z/3 \oplus \Z/3 &(0,0)& 2q^{-1}t^{-1} + 5 + 2qt \\
&&\pm (0,1)& q^{-2/3} (q^{-1}t^{-1} + 3 + qt) \\
&&\pm(1,0)& 1 \\
&&\pm(1,1)& 1 \\
&&\pm(1,2)& q^{-4/3} \\
\\

9_{47} &\Z/3 \oplus \Z/9 &(0,0)& q^{-3}t^{-3} + 4q^{-2}t^{-2} + 6q^{-1}t^{-1} + 5 + 6qt + 4q^2t^2 + q^3t^3 \\
&&\pm(0,1)& q^{-1/9} (q^{-1}t^{-1} + 3 + qt) \\
&&\pm(0,2)& q^{-4/9} (q^{-1}t^{-1} + 1 + qt) \\
&&\pm(0,3)& q^{-1}t^{-1} + 1 + qt \\
&&\pm(0,4)& q^{-7/9} \\
&&\pm(1,0)& q^{-1/3} \\
&&\pm(1,1)& q^{-1/9} (q^{-1}t^{-1} + 3 + qt) \\
&&\pm(1,2)& q^{-1/9} (q^{-1}t^{-1} + 3 + qt) \\
&&\pm(1,3)& q^{-1/3} \\
&&\pm(1,4)& q^{-7/9} \\
&&\pm(1,5)& q^{-4/9} (q^{-1}t^{-1} + 1 + qt) \\
&&\pm(1,6)& q^{-1/3} \\
&&\pm(1,7)& q^{-4/9} (q^{-1}t^{-1} + 1 + qt) \\
&&\pm(1,8)& q^{-7/9} \\
\end{array}
\]

\[
\begin{array}{cccl}
K & H_1(\Sigma^2(K);\Z) & \s & \sum_{i,j} \dim_{\Z/2}
\HFK_j(\Sigma^2(K), \tilde K, \s, i; \Z/2) t^i q^j \\
\\

9_{48} &\Z/3 \oplus \Z/9 &(0,0)& q^{-2}t^{-2} + 7q^{-1}t^{-1} + 11 + 7qt + q^2t^2 \\
&&\pm(0,1)& q^{-4/9} (q^{-1}t^{-1} + 1 + qt) \\
&&\pm(0,2)& q^{2/9} (2q^{-1}t^{-1} + 3 + 2qt) \\
&&\pm(0,3)& q^{-1}t^{-1} + 1 + qt \\
&&\pm(0,4)& q^{-1/9} \\
&&\pm(1,0)& q^{1/3} \\
&&\pm(1,1)& q^{2/9} (2q^{-1}t^{-1} + 3 + 2qt) \\
&&\pm(1,2)& q^{2/9} (2q^{-1}t^{-1} + 3 + 2qt) \\
&&\pm(1,3)& q^{1/3} \\
&&\pm(1,4)& q^{-4/9} (q^{-1}t^{-1} + 1 + qt) \\
&&\pm(1,5)& q^{-1/9} \\
&&\pm(1,6)& q^{1/3} \\
&&\pm(1,7)& q^{-1/9} \\
&&\pm(1,8)& q^{-4/9} (q^{-1}t^{-1} + 1 + qt) \\
\\

9_{49} &\Z/5 \oplus \Z/5 &(0,0)& 3q^{-2}t^{-2} + 6q^{-1}t^{-1} + 7 + 6qt + 3q^2t^2 \\
&&\pm(0,1)& q^{-2/5} (q^{-2}t^{-2} + q^{-1}t^{-1} + 1 + qt + q^2t^2) \\
&&\pm(0,2)& q^{2/5} \\
&&\pm(1,0)& q^{-2/5} (q^{-2}t^{-2} + q^{-1}t^{-1} + 1 + qt + q^2t^2) \\
&&\pm(1,1)& q^{-1/5} (q^{-1}t^{-1} + 1 + qt) \\
&&\pm(1,2)& q^{1/5} (2q^{-1}t^{-1} + 3 + 2qt) \\
&&\pm(1,3)& q^{1/5} (2q^{-1}t^{-1} + 3 + 2qt) \\
&&\pm(1,4)& q^{-2/5} (q^{-2}t^{-2} + q^{-1}t^{-1} + 1 + qt + q^2t^2) \\
&&\pm(2,0)& q^{2/5} \\
&&\pm(2,1)& q^{1/5} (2q^{-1}t^{-1} + 3 + 2qt) \\
&&\pm(2,2)& q^{-1/5} (q^{-1}t^{-1} + 1 + qt) \\
&&\pm(2,3)& q^{2/5} \\
&&\pm(2,4)& q^{-1/5} (q^{-1}t^{-1} + 1 + qt) \\
\\

10_{124} &\{0\} &0& q^{-4}t^{-4} + q^{-3}t^{-3} + \underline{t^{-1}
+ q +
q^2t} + q^3t^3 + q^4t^4 \\
\\

10_{128} &\Z/11 &0& 2q^{-3}t^{-3} + 3q^{-2}t^{-2} + q^{-1}t^{-1} +
\underline{q} + qt + 3q^2t^2 + 2q^3t^3 \\
&&\pm1& q^{8/11} (2q^{-1}t^{-1} + 3 + 2qt) \\
&&\pm2& q^{10/11} (q^{-1}t^{-1} + 1 + qt) \\
&&\pm3& q^{6/11} (q^{-1}t^{-1} + 1 + qt) \\
&&\pm4& q^{-4/11} (q^{-2}t^{-2} + q^{-1}t^{-1} + \underline{q} + qt + q^2t^2) \\
&&\pm5& q^{2/11} (q^{-1}t^{-1} + 1 + qt) \\
\\

10_{129} & \Z/25 &0& 2q^{-2}t^{-2} + 6q^{-1}t^{-1} + 9 + 6qt + 2q^2t^2 \\
&&\pm1& q^{-8/25} (q^{-2}t^{-2} + 2q^{-1}t^{-1} + 3 + 2qt + q^2t^2) \\
&&\pm2& q^{-7/25} (q^{-1}t^{-1} + 1 + qt) \\
&&\pm3& q^{3/25} (2q^{-1}t^{-1} + 3 + 2qt) \\
&&\pm4& q^{-3/25} (q^{-1}t^{-1} + 1 + qt) \\
&&\pm5& 1 \\
&&\pm6& q^{12/25} \\
&&\pm7& q^{8/25} \\
&&\pm8& q^{-12/25} \\
&&\pm9& q^{2/25} (q^{-1}t^{-1} + 3 + qt) \\
&&\pm10& 1 \\
&&\pm11& q^{7/25} (q^{-1}t^{-1} + 1 + qt) \\
&&\pm12& q^{23/25} (q^{-1}t^{-1} + 1 + qt) \\
\\

10_{130} & \Z/17 &0& 2q^{-2}t^{-2} + 4q^{-1}t^{-1} + 5 + 4qt + 2q^2t^2 \\
&&\pm1& q^{4/17} (q^{-2}t^{-2} + 2q^{-1}t^{-1} + 3 + 2qt + q^2t^2) \\
&&\pm2& q^{16/17} \\
&&\pm3& q^{19/17} (q^{-1}t^{-1} + 1 + qt) \\
&&\pm4& q^{13/17} (2q^{-1}t^{-1} + 3 + 2qt) \\
&&\pm5& q^{15/17} (q^{-1}t^{-1} + 1 + qt) \\
&&\pm6& q^{8/17} \\
&&\pm7& q^{9/17} (q^{-1}t^{-1} + 1 + qt) \\
&&\pm8& q^{1/17} (q^{-1}t^{-1} + 1 + qt) \\
\end{array}
\]

\[
\begin{array}{cccl}
K & H_1(\Sigma^2(K);\Z) & \s & \sum_{i,j} \dim_{\Z/2}
\HFK_j(\Sigma^2(K), \tilde K, \s, i; \Z/2) t^i q^j \\
\\

10_{131} &\Z/31 &0& 2q^{-2}t^{-2} + 8q^{-1}t^{-1} + 11 + 8qt + 2q^2t^2 \\
&&\pm1& q^{-18/31} (q^{-1}t^{-1} + 1 + qt) \\
&&\pm2& q^{-10/31} (q^{-1}t^{-1} + 1 + qt) \\
&&\pm3& q^{-7/31} (q^{-1}t^{-1} + 3 + qt) \\
&&\pm4& q^{-9/31} \\
&&\pm5& q^{15/31} \\
&&\pm6& q^{3/31} \\
&&\pm7& q^{-14/31} (q^{-1}t^{-1} + 1 + qt) \\
&&\pm8& q^{-5/31} (2q^{-1}t^{-1} + 5 + 2qt) \\
&&\pm9& q^{-1/31} \\
&&\pm10& q^{-2/31} (q^{-1}t^{-1} + 1 + qt) \\
&&\pm11& q^{-8/31} (q^{-2}t^{-2} + 4q^{-1}t^{-1} + 5 + 4qt + q^2t^2) \\
&&\pm12& q^{-19/31} (q^{-1}t^{-1} + 3 + qt) \\
&&\pm13& q^{-4/31} (2q^{-1}t^{-1} + 3 + 2qt) \\
&&\pm14& q^{-25/31} \\
&&\pm15& q^{11/31} \\
\\

10_{132} &\Z/5 &0& q^{-2}t^{-2} + (2q^{-1}+ \underline{1})t^{-1} +
(2+\underline{q}) + (2q+\underline{q^2})t + q^2t^2 \\
&&\pm1& q^{2/5} \\
&&\pm2& q^{3/5} (q^{-1}t^{-1} + 1 + qt) \\
\\

10_{133} &\Z/19 &0& q^{-2}t^{-2} + 5q^{-1}t^{-1} + 7 + 5qt + q^2t^2 \\
&&\pm1& q^{-3/19} \\
&&\pm2& q^{-12/19} (q^{-1}t^{-1} + 1 + qt) \\
&&\pm3& q^{-8/19} (q^{-1}t^{-1} + 1 + qt) \\
&&\pm4& q^{9/19} \\
&&\pm5& q^{1/19} \\
&&\pm6& q^{-13/19} (q^{-1}t^{-1} + 3 + qt) \\
&&\pm7& q^{5/19} \\
&&\pm8& q^{-2/19} (2q^{-1}t^{-1} + 3 + 2qt) \\
&&\pm9& q^{-15/19} \\
\\

10_{134} &\Z/23 &0& 2q^{-3}t^{-3} + 4q^{-2}t^{-2} + 4q^{-1}t^{-1} +
3 + 4qt + 4q^2t^2 + 2q^3t^3 \\
&&\pm1& q^{8/23} (q^{-1}t^{-1} + 1 + qt) \\
&&\pm2& q^{9/23} (q^{-2}t^{-2} + q^{-1}t^{-1} + 1 + qt + q^2t^2) \\
&&\pm3& q^{3/23} (q^{-2}t^{-2} + 2q^{-1}t^{-1} + 3 + 2qt + q^2t^2) \\
&&\pm4& q^{-10/23} (q^{-3}t^{-3} + q^{-2}t^{-2} + \underline{q} + q^2t^2 + q^3t^3) \\
&&\pm5& q^{16/23} (q^{-1}t^{-1} + 1 + qt) \\
&&\pm6& q^{12/23} (q^{-1}t^{-1} + 1 + qt) \\
&&\pm7& q^{1/23} (q^{-2}t^{-2} + q^{-1}t^{-1} + 1 + qt + q^2t^2) \\
&&\pm8& q^{29/23} \\
&&\pm9& q^{4/23} (q^{-1}t^{-1} + 1 + qt) \\
&&\pm10& q^{18/23} (2q^{-1}t^{-1} + 3 + 2qt) \\
&&\pm11& q^{25/23} \\
\\

10_{135} &\Z/37 &0& 3q^{-2}t^{-2} + 9q^{-1}t^{-1} + 13 + 9qt +
3q^2t^2 \\
&&\pm1& q^{14/37} \\
&&\pm2& q^{-18/37} \\
&&\pm3& q^{15/37} (2q^{-1}t^{-1} + 3 + 2qt) \\
&&\pm4& q^{2/37} \\
&&\pm5& q^{17/37} (q^{-1}t^{-1} + 1 + qt) \\
&&\pm6& q^{-14/37} (q^{-2}t^{-2} + 2q^{-1}t^{-1} + 3 + 2qt + q^2t^2) \\
&&\pm7& q^{-17/37} (q^{-1}t^{-1} + 1 + qt) \\
&&\pm8& q^{8/37} (q^{-2}t^{-2} + 3q^{-1}t^{-1} + 5 + 3qt + q^2t^2) \\
&&\pm9& q^{-13/37} (q^{-1}t^{-1} + 1 + qt) \\
&&\pm10& q^{-6/37} \\
&&\pm11& q^{29/37} (q^{-1}t^{-1} + 1 + qt) \\
&&\pm12& q^{18/37} \\
&&\pm13& q^{-2/37} \\
&&\pm14& q^{6/37} (q^{-2}t^{-2} + 2q^{-1}t^{-1} + 3 + 2qt + q^2t^2) \\
&&\pm15& q^{5/37} (2q^{-1}t^{-1} + 3 + 2qt) \\
&&\pm16& q^{-5/37} (q^{-1}t^{-1} + 1 + qt) \\
&&\pm17& q^{13/37} (q^{-1}t^{-1} + 1 + qt) \\
&&\pm18& q^{22/37} \\
\\
\end{array}
\]

\[
\begin{array}{cccl}
K & H_1(\Sigma^2(K);\Z) & \s & \sum_{i,j} \dim_{\Z/2}
\HFK_j(\Sigma^2(K), \tilde K, \s, i; \Z/2) t^i q^j \\
\\

10_{136} &\Z/15 &1& q^{-2}t^{-2} + 4q^{-1}t^{-1} + 6 + \underline{q} + 4qt + q^2t^2 \\
&&\pm1& q^{7/15} \\
&&\pm2& q^{13/15} (q^{-1}t^{-1} + 3 + qt) \\
&&\pm3& q^{1/5} \\
&&\pm4& q^{7/15} \\
&&\pm5& q^{2/3} (q^{-1}t^{-1} + 1 + qt) \\
&&\pm6& q^{4/5} (2q^{-1}t^{-1} + 3 + 2qt) \\
&&\pm7& q^{13/15} (q^{-1}t^{-1} + 3 + qt) \\
\\

10_{139} &\Z/3 &0& q^{-4}t^{-4} + q^{-3}t^{-3} + \underline{2qt^{-1}
+ 3q + 2q^3t} + q^3t^3 + q^4t^4 \\
&&\pm1& q^{5/3}( q^{-2}t^{-2} + q^{-1}t^{-1} + 1 + qt + q^2t^2 ) \\
\\

10_{140} &\Z/9 &0& q^{-2}t^{-2} + 2q^{-1}t^{-1} + 3 + 2qt + q^2t^2 \\
&&\pm1& q^{11/9} (q^{-1}t^{-1} + 1 + qt) \\
&&\pm2& q^{8/9} \\
&&\pm3& 1 \\
&&\pm4& q^{5/9} (q^{-1}t^{-1} + 1 + qt) \\
\\

10_{142} &\Z/15 &0& q^{-3}t^{-3} + 3q^{-2}t^{-2} + 2q^{-1}t^{-1} +
1 + 2qt + 3q^2t^2 + 2q^3t^3 \\
&&\pm1& q^{1/15} (q^{-2}t^{-2} + q^{-1}t^{-1} + 1 + qt + q^2t^2) \\
&&\pm2& q^{4/15} (q^{-1}t^{-1} + 1 + qt) \\
&&\pm3& q^{-2/5} (q^{-3}t^{-3} + q^{-2}t^{-2} + \underline{q} + q^2t^2 + q^3t^3) \\
&&\pm1& q^{1/15} (q^{-2}t^{-2} + q^{-1}t^{-1} + 1 + qt + q^2t^2) \\
&&\pm6& q^{2/3} (2q^{-1}t^{-1} + 3 + 2qt) \\
&&\pm2& q^{7/5} \\
&&\pm2& q^{4/15} (q^{-1}t^{-1} + 1 + qt) \\
\\

10_{145} &\Z/3 &0& q^{-2}t^{-2} + (q^{-1}+\underline{2q})t^{-1} +
\underline{q+4q^2} + (q+\underline{2q^3})t + q^2t^2 \\
&&\pm1& q^{4/3} (2q^{-1}t^{-1} + 3 + 2qt) \\
\\

10_{147} &\Z/27 &0& 2q^{-2}t^{-2} + 7q^{-1}t^{-1} + 9 + 7qt +
2q^2t^2 \\
&&\pm1& q^{7/27} (q^{-1}t^{-1} + 3 + qt) \\
&&\pm2& q^{1/27} \\
&&\pm3& q^{1/3} (2q^{-1}t^{-1} + 5 + 2qt) \\
&&\pm4& q^{4/27} (q^{-2}t^{-2} + 3q^{-1}t^{-1} + 3 + 3qt + q^2t^2) \\
&&\pm5& q^{13/27} \\
&&\pm6& q^{1/3} \\
&&\pm7& q^{19/27} (q^{-1}t^{-1} + 3 + qt) \\
&&\pm8& q^{16/27} (q^{-1}t^{-1} + 1 + qt) \\
&&\pm9& q^{-1}t^{-1} + 1 + qt \\
&&\pm10& q^{25/27} \\
&&\pm11& q^{37/27} \\
&&\pm12& q^{1/3} \\
&&\pm13& q^{22/27} (2q^{-1}t^{-1} + 3 + 2qt) \\
\end{array}
\]

\[
\begin{array}{cccl}
K & H_1(\Sigma^2(K);\Z) & \s & \sum_{i,j} \dim_{\Z/2}
\HFK_j(\Sigma^2(K), \tilde K, \s, i; \Z/2) t^i q^j \\
\\

10_{158} &\Z/45 &0& q^{-3}t^{-3} + 4q^{-2}t^{-2} + 10q^{-1}t^{-1} +
15 + 10qt + 4q^2t^2 + q^3t^3 \\
&&\pm1& q^{8/45} (q^{-1}t^{-1} + 3 + qt) \\
&&\pm2& q^{-13/45} (q^{-2}t^{-2} + 3q^{-1}t^{-1} + 3 + 3qt + q^2t^2) \\
&&\pm3& q^{-2/5} \\
&&\pm4& q^{38/45} \\
&&\pm5& q^{4/9} (q^{-1}t^{-1} + 3 + qt) \\
&&\pm6& q^{-2/5} \\
&&\pm7& q^{-13/45} (q^{-2}t^{-2} + 3q^{-1}t^{-1} + 3 + 3qt + q^2t^2) \\
&&\pm8& q^{17/45} (q^{-2}t^{-2} + 3q^{-1}t^{-1} + 3 + 3qt + q^2t^2) \\
&&\pm9& q^{2/5} (2q^{-1}t^{-1} + 5 + 2qt) \\
&&\pm10& q^{-2/9} (2q^{-1}t^{-1} + 5 + 2qt) \\
&&\pm11& q^{-22/45} \\
&&\pm12& q^{-2/5} \\
&&\pm13& q^{-2/45} \\
&&\pm14& q^{38/45} \\
&&\pm15& q^{-1}t^{-1} + 3 + qt \\
&&\pm16& q^{-22/45} \\
&&\pm17& q^{17/45} (q^{-2}t^{-2} + 3q^{-1}t^{-1} + 3 + 3qt + q^2t^2) \\
&&\pm18& q^{-2/5} \\
&&\pm19& q^{8/45} (q^{-1}t^{-1} + 3 + qt) \\
&&\pm20& q^{1/9} (q^{-1}t^{-1} + 1 + qt) \\
&&\pm21& q^{2/5} \\
&&\pm22& q^{2/45} \\
\\

10_{160} &\Z/21 &0&  q^{-3}t^{-3} + 4q^{-2}t^{-2} + 4q^{-1}t^{-1} +
3 + 4qt + 4q^2t^2 + q^3t^3 \\
&&\pm1& q^{1/21} (q^{-1}t^{-1} + 1 + qt) \\
&&\pm2& q^{4/21} (q^{-2}t^{-2} + q^{-1}t^{-1} + 1 + qt + q^2t^2) \\
&&\pm3& q^{3/7} (q^{-1}t^{-1} + 1 + qt) \\
&&\pm4& q^{16/21} \\
&&\pm5& q^{4/21} (q^{-2}t^{-2} + q^{-1}t^{-1} + 1 + qt + q^2t^2) \\
&&\pm6& q^{5/7} (2q^{-1}t^{-1} + 3 + 2qt) \\
&&\pm7& q^{4/3} \\
&&\pm8& q^{1/21} (q^{-1}t^{-1} + 1 + qt) \\
&&\pm9& q^{6/7} (q^{-1}t^{-1} + 3 + qt) \\
&&\pm10& q^{16/21} \\
\\

10_{161} &\Z/5 &0& q^{-3}t^{-3} + (q^{-2}+\underline{1})t^{-2} +
\underline{2qt^{-1} + 3q^2 + 2q^3t} + (q^2+\underline{q^4})t^2 + q^3t^3 \\
&&\pm1& q^{9/5} (2q^{-1}t^{-1} + 3 + 2qt) \\
&&\pm2& q^{6/5} (q^{-2}t^{-2} + q^{-1}t^{-1} + 1 + qt + q^2t^2) \\
\\

10_{164} &\Z/45 &0& 3q^{-2}t^{-2} + 11q^{-1}t^{-1} + 17 + 11qt +
3q^2t^2 \\
&&\pm1& q^{17/45} (q^{-1}t^{-1} + 1 + qt) \\
&&\pm2& q^{-22/45} \\
&&\pm3& q^{2/5} \\
&&\pm4& q^{2/45} \\
&&\pm5& q^{4/9} (q^{-1}t^{-1} + 3 + qt) \\
&&\pm6& q^{-2/5} \\
&&\pm7& q^{-2/45} \\
&&\pm8& q^{8/45} (q^{-2}t^{-2} + 3q^{-1}t^{-1} + 5 + 3qt + q^2t^2) \\
&&\pm9& q^{-2/5} (q^{-2}t^{-2} + 2q^{-1}t^{-1} + 3 + 2qt + q^2t^2) \\
&&\pm10& q^{-2/9} \\
&&\pm11& q^{-13/45} (q^{-1}t^{-1} + 1 + qt) \\
&&\pm12& q^{2/5} \\
&&\pm13& q^{38/45} \\
&&\pm14& q^{2/45} \\
&&\pm15& q^{-1}t^{-1} + 3 + qt \\
&&\pm16& q^{-13/45} (q^{-1}t^{-1} + 1 + qt) \\
&&\pm17& q^{8/45} (q^{-2}t^{-2} + 3q^{-1}t^{-1} + 5 + 3qt + q^2t^2) \\
&&\pm18& q^{2/5} \\
&&\pm19& q^{17/45} (q^{-1}t^{-1} + 1 + qt) \\
&&\pm20& q^{1/9} (2q^{-1}t^{-1} + 3 + 2qt) \\
&&\pm21& q^{-2/5} \\
&&\pm22& q^{38/45} \\
\end{array}
\]

\[
\begin{array}{cccl}
K & H_1(\Sigma^2(K);\Z) & \s & \sum_{i,j} \dim_{\Z/2}
\HFK_j(\Sigma^2(K), \tilde K, \s, i; \Z/2) t^i q^j \\
\\

11n_{12} &\Z/13 &0& q^{-2}t^{-2} +
(\underline{q^{-2}}+4q^{-1})t^{-1} + \underline{q^{-1}}+6 +
(\underline{1}+4q)t + q^2t^2 \\
&&\pm1& q^{-2/13} \\
&&\pm2& q^{-8/13} (q^{-1}t^{-1} + 3 + qt) \\
&&\pm3& q^{-18/13} \\
&&\pm4& q^{-6/13} \\
&&\pm5& q^{-11/13} (2q^{-1}t^{-1} + 3 + 2qt) \\
&&\pm6& q^{-7/13} (q^{-1}t^{-1} + 1 + qt) \\
\\

11n_{19} &\Z/5 &0& q^{-3}t^{-3} + 2q^{-2}t^{-2} +
(q^{-1}+\underline{1})
t^{-1} + \underline{q} + (q+\underline{q^2}) t + 2q^2t^2 + q^3t^3 \\
&&\pm1& q^{4/5} (q^{-2}t^{-2} + q^{-1}t^{-1} + 1 + qt + q^2t^2) \\
&&\pm2& q^{6/5}(q^{-1}t^{-1} + 3 + qt) \\
\\

11n_{20} &\Z/23 &0& 2q^{-2}t^{-2} + 6q^{-1}t^{-1} +
\underline{8+q} + 6qt + 2q^2t^2 \\
&&\pm1& q^{17/23} (q^{-1}t^{-1} + 3 + qt) \\
&&\pm2& q^{-1/23} \\
&&\pm3& q^{-8/23} q^{-2}t^{-2} + 2q^{-1}t^{-1} + \underline{2+q} + 2qt + 1q^2t^2 \\
&&\pm4& q^{19/23} (2q^{-1}t^{-1} + 5 + 2qt) \\
&&\pm5& q^{11/23} \\
&&\pm6& q^{14/23} (q^{-1}t^{-1} + 1 + qt) \\
&&\pm7& q^{5/23} (q^{-1}t^{-1} + 3 + qt) \\
&&\pm8& q^{7/23} \\
&&\pm9& q^{20/23} (2q^{-1}t^{-1} + 3 + 2qt) \\
&&\pm10& q^{21/23} (q^{-1}t^{-1} + 3 + qt) \\
&&\pm11& q^{10/23} (q^{-1}t^{-1} + 1 + qt) \\
\\

11n_{38} &\Z/3 &0& q^{-2}t^{-2} + (2q^{-1}+\underline{1})t^{-1} +
\underline{2+3q} + (2q+\underline{q^2})t + q^2t^2 \\
&&\pm1& q^{4/3} (q^{-1}t^{-1} + 1 + qt) \\
\\

11n_{49} &\{0\} &0& q^{-2}t^{-2} +
(\underline{4q^{-3}}+2q^{-1})t^{-1}
+ \underline{9q^{-2}}+2 + (\underline{4q^{-1}}+q)t + q^2t^2 \\
\\

11n_{95} &\Z/33 &0& q^{-3}t^{-3} + 5q^{-2}t^{-2} + 7q^{-1}t^{-1} + 7
+ 7qt + 5q^2t^2 + q^3t^3 \\
&&\pm1& q^{-13/33} (q^{-1}t^{-1} + 1 + qt) \\
&&\pm2& q^{14/33} \\
&&\pm3& q^{5/11} (q^{-1}t^{-1} + 1 + qt) \\
&&\pm4& q^{-10/33} (q^{-2}t^{-2} + q^{-1}t^{-1} + 1 + qt + q^2t^2) \\
&&\pm5& q^{5/33} (2q^{-1}t^{-1} + 3 + 2qt) \\
&&\pm6& q^{-2/11} (q^{-2}t^{-2} + q^{-1}t^{-1} + 1 + qt + q^2t^2) \\
&&\pm7& q^{-10/33} (q^{-2}t^{-2} + q^{-1}t^{-1} + 1 + qt + q^2t^2) \\
&&\pm8& q^{26/33} \\
&&\pm9& q^{1/11} (q^{-1}t^{-1} + 1 + qt) \\
&&\pm10& q^{-13/33} (q^{-1}t^{-1} + 1 + qt) \\
&&\pm11& q^{1/3} (q^{-1}t^{-1} + 1 + qt) \\
&&\pm12& q^{3/11} (2q^{-1}t^{-1} + 3 + 2qt) \\
&&\pm13& q^{14/33} \\
&&\pm14& q^{26/33} \\
&&\pm15& q^{4/11} (q^{-1}t^{-1} + 3 + qt) \\
&&\pm16& q^{5/33} (2q^{-1}t^{-1} + 3 + 2qt) \\
\\

11n_{102} &\Z/3 &0& q^{-2}t^{-2} + (5q^{-1}+\underline{2q})t^{-1}
+ 7+ \underline{4q^{2}} + (5q+\underline{2q^{3}})t + q^2t^2 \\
&&\pm1& q^{1/3} (2q^{-1}t^{-1}+5+2qt) \\
\\

11n_{116} &\{0\} &0& q^{-2}t^{-2} +
(\underline{4q^{-3}}+2q^{-1})t^{-1}
+ \underline{9q^{-2}}+2 + (\underline{4q^{-1}}+q)t + q^2t^2 \\
\end{array}
\]

\[
\begin{array}{cccl}
K & H_1(\Sigma^2(K);\Z) & \s & \sum_{i,j} \dim_{\Z/2}
\HFK_j(\Sigma^2(K), \tilde K, \s, i; \Z/2) t^i q^j \\
\\

11n_{117} &\Z/35 &0& 3q^{-2}t^{-2} + 9q^{-1}t^{-1}+11 + 9qt +
3q^2t^2 \\
&&\pm1& q^{9/35} (2q^{-1}t^{-1} + 5 + 2qt) \\
&&\pm2& q^{1/35} \\
&&\pm3& q^{11/35} (q^{-1}t^{-1} + 3 + qt) \\
&&\pm4& q^{4/35} (q^{-2}t^{-2} + 3q^{-1}t^{-1} + 3 + 3qt + q^2t^2) \\
&&\pm5& q^{3/7} (q^{-1}t^{-1} + 3 + qt) \\
&&\pm6& q^{9/35} (2q^{-1}t^{-1} + 5 + 2qt) \\
&&\pm7& q^{-2/5} (q^{-2}t^{-2} + 2q^{-1}t^{-1} + 2 + \underline{q} + 2qt + q^2t^2) \\
&&\pm8& q^{16/35} (q^{-1}t^{-1} + 1 + qt) \\
&&\pm9& q^{29/35} \\
&&\pm10& q^{5/7} (2q^{-1}t^{-1} + 5 + 2qt) \\
&&\pm11& q^{4/35} (q^{-2}t^{-2} + 3q^{-1}t^{-1} + 3 + 3qt + q^2t^2) \\
&&\pm12& q^{1/355} \\
&&\pm13& q^{16/35} (q^{-1}t^{-1} + 1 + qt) \\
&&\pm14& q^{7/5} \\
&&\pm15& q^{6/7} (2q^{-1}t^{-1} + 3 + 2qt) \\
&&\pm16& q^{29/35} \\
&&\pm17& q^{11/35} (q^{-1}t^{-1} + 3 + qt) \\
\\

11n_{118} &\Z/21 &0& q^{-3}t^{-3} + 4q^{-2}t^{-2} + 4q^{-1}t^{-1} +
3 + 4qt + 4q^2t^2 + q^3t^3 \\
&&\pm1& q^{5/21} (q^{-1}t^{-1} + 1 + qt) \\
&&\pm2& q^{20/21} \\
&&\pm3& q^{1/7} (2q^{-1}t^{-1} + 3 + 2qt) \\
&&\pm4& q^{-4/21} (q^{-2}t^{-2} + q^{-1}t^{-1} + 1 + qt + q^2t^2) \\
&&\pm5& q^{20/21} \\
&&\pm6& q^{4/7} \\
&&\pm7& q^{-1/3} (q^{-1}t^{-1} + 1 + qt) \\
&&\pm8& q^{5/21} (q^{-1}t^{-1} + 1 + qt) \\
&&\pm9& q^{2/7} (q^{-1}t^{-1} + 3 + qt) \\
&&\pm10& q^{-4/21} (q^{-2}t^{-2} + q^{-1}t^{-1} + 1 + qt + q^2t^2) \\
\\

11n_{122} &\Z/27 &0& 2q^{-2}t^{-2} + 7q^{-1}t^{-1} + 9 + 2qt +
2q^2t^2 \\
&&\pm1& q^{13/27} \\
&&\pm2& q^{-2/27} (2q^{-1}t^{-1} + 3 + 2qt) \\
&&\pm3& q^{1/3} (2q^{-1}t^{-1} + 5 + 2qt) \\
&&\pm4& q^{-8/27} (q^{-1}t^{-1} + 1 + qt) \\
&&\pm5& q^{1/27} \\
&&\pm6& q^{1/3} \\
&&\pm7& q^{-11/27} \\
&&\pm8& q^{-5/27} (q^{-1}t^{-1} + 3 + qt) \\
&&\pm9& q^{-1}t^{-1} + 1 + qt \\
&&\pm10& q^{-23/27} \\
&&\pm11& q^{-20/27} (q^{-2}t^{-2} + 3q^{-1}t^{-1} + 3 + 3qt + q^2t^2) \\
&&\pm12& q^{1/3} \\
&&\pm13& q^{-17/27} (q^{-1}t^{-1} + 3 + qt) \\
\\

11n_{138} &\Z/15 &0& 2q^{-2}t^{-2} + 4q^{-1}t^{-1} +
(\underline{q^{-1}} + 4) + 4qt + 2q^2t^2 \\
&&\pm1& q^{-7/15} \\
&&\pm2& q^{-13/15} (q^{-1}t^{-1} + 3 + qt) \\
&&\pm3& q^{-1/5} ((q^{-2} + 2q^{-1})t^{-1} + (q^{-1} + 4) + (1+2q)t
) \\
&&\pm4& q^{-7/15} \\
&&\pm5& q^{-2/3} (q^{-2}t^{-2} + 3q^{-1}t^{-1} + 3 + 3qt + q^2t^2) \\
&&\pm6& q^{-9/5} \\
&&\pm7& q^{-13/15} (q^{-1}t^{-1} + 3 + qt) \\
\\

11n_{139} &\Z/9 &0& 2q^{-1}t^{-1} + 5 + 2qt \\
&&\pm 1& q^{-4/9} \\
&&\pm 2& q^{-16/9} \\
&&\pm 3& 1 \\
&&\pm 4& q^{-10/9} (q^{-1}t^{-1} + 3 + qt) \\
\end{array}
\]

\[
\begin{array}{cccl}
K & H_1(\Sigma^2(K);\Z) & \s & \sum_{i,j} \dim_{\Z/2}
\HFK_j(\Sigma^2(K), \tilde K, \s, i; \Z/2) t^i q^j \\
\\

11n_{141} &\Z/21 &0& 5q^{-1}t^{-1} + 11 + 5qt \\
&&\pm1& q^{-10/21} \\
&&\pm2& q^{2/21} (2q^{-1}t^{-1} + 5 + 2qt) \\
&&\pm3& q^{-2/7} (2q^{-1}t^{-1} + 5 + 2qt) \\
&&\pm4& q^{8/21} (q^{-1}t^{-1} + 3 + qt) \\
&&\pm5& q^{2/21} (2q^{-1}t^{-1} + 5 + 2qt) \\
&&\pm6& q^{6/7} \\
&&\pm7& q^{2/3} (2q^{-1}t^{-1} + 5 + 2qt) \\
&&\pm8& q^{-10/21} \\
&&\pm9& q^{10/7} \\
&&\pm10& q^{8/21} (q^{-1}t^{-1} + 3 + qt) \\
\\

11n_{142} &\Z/33 &0& q^{-2}t^{-2} + 8q^{-1}t^{-1} +
15 + 8qt + q^2t^2 \\
&&\pm1& q^{2/33} (q^{-1}t^{-1} + 3 + qt) \\
&&\pm2& q^{8/33} (2q^{-1}t^{-1} + 5 + 2qt) \\
&&\pm3& q^{6/11} (q^{-1}t^{-1} + 3 + qt) \\
&&\pm4& q^{32/33} \\
&&\pm5& q^{-16/33} \\
&&\pm6& q^{2/11} (q^{-1}t^{-1} + 3 + qt) \\
&&\pm7& q^{32/33} \\
&&\pm8& q^{-4/33} \\
&&\pm9& q^{10/11} (q^{-1}t^{-1} + 3 + qt) \\
&&\pm10& q^{2/33} (q^{-1}t^{-1} + 3 + qt) \\
&&\pm11& q^{4/3} \\
&&\pm12& q^{8/11} (2q^{-1}t^{-1} + 5 + 2qt) \\
&&\pm13& q^{8/33} (2q^{-1}t^{-1} + 5 + 2qt) \\
&&\pm14& q^{-4/33} \\
&&\pm15& q^{-4/11} (2q^{-1}t^{-1} + 5 + 2qt) \\
&&\pm16& q^{-16/33} \\
\\

11n_{143} &\Z/9 &0& q^{-3}t^{-3} + (\underline{q^{-4}}+
3q^{-2})t^{-2} + (\underline{2q^{-3}}+3q^{-1})t^{-1} +
(\underline{2q^{-2}}+3) \\
&&& \quad + (\underline{2q^{-1}}+3q)t + (\underline{1}
+ 3q^2)t^2 + q^3t^3 \\
&&\pm1& q^{-10/9} ((q^{-1}+1)t^{-1} + (2+q) + (1+q^2)t) \\
&&\pm2& q^{-4/9} \\
&&\pm3& 1 \\
&&\pm4& q^{-7/9} (q^{-2}t^{-2} + 3q^{-1}t^{-1} + 3 + 3qt + q^2t^2) \\
\\

11n_{145} &\Z/9 &0& q^{-3}t^{-3} + (2q^{-2}+ \underline{1})t^{-2} +
(q^{-1}+\underline{4q})t^{-1} + \underline{7q^2} + (q+\underline{4q^3})t \\
&&& \quad + (2q^2+\underline{t^4})t^2 + q^3t^3 \\
&&\pm1& q^{10/9} (q^{-2}t^{-2} + 3q^{-1}t^{-1} + 5 + 3qt + q^2t^2) \\
&&\pm2& q^{22/9} \\
&&\pm3& q^2 \\
&&\pm4& q^{16/9} (q^{-2}t^{-2} + 3q^{-1}t^{-1} + 5 + 3qt + q^2t^2) \\

\end{array}
\]

\end{tiny}

\section{Observations}

Grigsby \cite{G1} showed that when $K \subset S^3$ is a two-bridge
knot, the Heegaard Floer knot homology of $\tilde K \subset
\Sigma^2(K)$ in the canonical spin$^c$ structure is isomorphic as a
bigraded $\Z/2$-vector space to that of $K \subset S^3$: i.e.,
$\HFK(\Sigma^2(K), \tilde K, \s_0) \cong \HFK(S^3,K)$, up to an
overall shift in the Maslov grading. Our results suggest that the
same is true for a wider class for knots. Specifically, we say that
$\HFK(S^3,K)$ is \emph{perfect} if it is supported along a single
diagonal, i.e., there exists a constant $C$ such that
$\HFK_j(S^3,K,i)=0$ when $j-i\ne C$. We conjecture:
\begin{conj} \label{conj:perfect}
Let $K \subset S^3$ be a knot such that $\HFK(S^3,K)$ is supported
along a single diagonal, i.e., Then $\HFK(\Sigma^2(K),\tilde K,
\s_0) \cong \HFK(S^3,K)$ as bigraded groups, up to a possible shift
in the absolute Maslov grading.
\end{conj}

It is well-known \cite{OSzAlt, R2} that $\HFK(S^3,K)$ is perfect
whenever $K$ is alternating (and hence for all two-bridge knots).
More generally, let $\QQ$ be the smallest set of link types such
that:
\begin{itemize}
\item The unknot is in $\QQ$.
\item Suppose $L$ admits a projection such that the two resolutions
at some crossing, $L_0$ and $L_1$, are both in $\QQ$ and satisfy
$\det(L_0) + \det(L_1) = \det(L)$. Then $L$ is in $\QQ$.
\end{itemize}
The links in $\QQ$ are called \emph{quasi-alternating}; for
instance, any alternating link is quasi-alternating. Manolescu and
Ozsv\'ath \cite{MO} have shown that whenever $L$ is
quasi-alternating, both $\HFK(S^3,L)$ and the Khovanov homology of
$L$ are perfect. (Additionally, Ozsv\'ath and Szab\'o
\cite{OSzDouble} have shown that the branched double cover of any
quasi-alternating link $L$ is an \emph{L-space}, meaning that
$\HF(\Sigma^2(L),\s)$ has rank 1 in each spin$^c$ structure.)
Conjecture \ref{conj:perfect} would then imply that
$\HFK(\Sigma^2(K),\tilde K,\s_0)$ is perfect whenever $K$ is
quasi-alternating.

One can also ask under what conditions $\HFK(\Sigma^2(K),\tilde
K,\s)$ is perfect when $\s \ne \s_0$. The knots $10_{134}$ and
$11n_{117}$ have the property that both $\HFK(S^3,K)$ and
$\HFK(\Sigma^2(K),\tilde K, \s_0)$ are perfect and isomorphic, but
there is a spin$^c$ structure $\s$ in which $\HFK(S^2(K),\tilde
K,\s)$ is not perfect. It is not known, however, whether these knots
are quasi-alternating.

On the other hand, when $\HFK(S^3,K)$ is not perfect, the
isomorphism between $\HFK(S^3,K)$ and $\HFK(\Sigma^2(K),\tilde K,
\s_0)$ fails. A few patterns are worth mentioning. If
$\HFK(S^3,K,g)$ (where $g=g(K)$) is supported in a single Maslov
grading $g+c$, define the \emph{main diagonal} of $\HFK(S^3,K)$ as
the groups $\HFK_{i+c}(S^3,K,i)$. (This assumption fails when the
rank of $\Delta_K$ is less than twice $g(K)$, for instance.) In
every example considered here, the remaining nonzero groups lie
either all above $(M>A+c)$ or all below $(M<A+c)$ the main diagonal.
(See \cite{BG} for the values of $\HFK(S^3,K)$ for all
non-alternating knots with $\le 12$ crossings.)

In most of our examples, the main diagonal of $\HFK(\Sigma^2(K),
\tilde K, \s_0)$ is isomorphic to that of $\HFK(S^3,K)$, while the
Maslov gradings of the off-diagonal groups may be shifted by an
overall constant. That constant is sometimes odd, implying that the
Maslov $\Z/2$-gradings need not be the same. For instance, when $K$
is the knot $10_{161}$, the off-diagonal groups are shifted by
three. However, there are also instances where the main diagonals of
$\HFK(S^3,K)$ and $\HFK(\Sigma^2(K), \tilde K, \s_0)$ are not
isomorphic. When $K = 10_{145}$, the matrices of the ranks of
$\HFK_j(S^3,K,i)$ and $\HFK_j(\Sigma^2(K), \tilde K,i))$ are,
respectively,
\[
\begin{pmatrix}
0 &     0 &     0 &     0 & \bm 1 \\
0 &     0 &     0 & \bm 1 &     0 \\
0 &     0 & \bm 1 &     2 &     0 \\
0 & \bm 1 &     4 &     0 &     0 \\
\bm 1 & 2 &     0 &     0 &     0 \\
\end{pmatrix}
\text{ and }
\begin{pmatrix}
0 &     0 &     0 &     2 &     0 \\
0 &     0 &     4 &     0 & \bm 1 \\
0 &     2 &     1 & \bm 1 &     0 \\
0 &     0 & \bm 0 &     0 &     0 \\
0 & \bm 1 &     0 &     0 &     0 \\
\bm 1 & 0 &     0 &     0 &     0 \\
\end{pmatrix}
\]
(where the Alexander grading is on the horizontal axis, the Maslov
grading is on the vertical axis, and the main diagonal is shown in
bold). Here, one of the groups on the main diagonal in $\HFK(S^3,K)$
is shifted upward by one. In this case, the total rank in each
Alexander grading is still the same, but there are also instances
where that statement fails to hold. For the knots $11n_{49}$ and
$11n_{116}$, which have determinant 1 and identical Heegaard Floer
homology both downstairs and upstairs, the ranks of
$\HFK_j(S^3,K,i)$ and $\HFK_j(\Sigma^2(K),\tilde K,\s_0,i)$ (in the
unique spin$^c$ structure) are given by
\[
\begin{pmatrix}
0 &     0 &     0 &     2 & \bm 1 \\
0 &     0 &     5 & \bm 2 &     0 \\
0 &     2 & \bm 2 &     0 &     0 \\
0 & \bm 2 &     0 &     0 &     0 \\
\bm 1 & 0 &     0 &     0 &     0 \\
\end{pmatrix}
\text{ and }
\begin{pmatrix}
0 &     0 &     0 &     0 & \bm 1 \\
0 &     0 &     0 & \bm 2 &     0 \\
0 &     0 & \bm 2 &     0 &     0 \\
0 & \bm 2 &     0 &     4 &     0 \\
\bm 1 & 0 &     9 &     0 &     0 \\
0 &     4 &     0 &     0 &     0 \\
\end{pmatrix}.
\]
Another example in which the total ranks of $\HFK(S^3,K)$ and
$\HFK(\Sigma^2(K),\tilde K, \s_0)$ are different is the knot
$11n_{102}$, for which the ranks are
\[
\begin{pmatrix}
0 &     0 &     0 &     2 & \bm 1 \\
0 &     0 &     4 & \bm 3 &     0 \\
0 &     2 & \bm 3 &     0 &     0 \\
0 & \bm 3 &     0 &     0 &     0 \\
\bm 1 & 0 &     0 &     0 &     0 \\
\end{pmatrix}
\text{ and }
\begin{pmatrix}
0 &     0 &     0 &     2 &     0 \\
0 &     0 &     4 &     0 & \bm 1 \\
0 &     2 &     0 & \bm 5 &     0 \\
0 &     0 & \bm 7 &     0 &     0 \\
0 & \bm 5 &     0 &     0 &     0 \\
\bm 1 & 0 &     0 &     0 &     0 \\
\end{pmatrix}.
\]

Finally, note that the pretzel knots $8_{20} = P(3,-3,2)$ and
$10_{140} = P(4,3,-3)$ have identical knot Floer homology but can be
distinguished by $\HFK(\Sigma^2(K), \tilde K)$. The relative Maslov
gradings between spin$^c$ structures are necessary in this case. For
another such example, see \cite{G1}.

\end{document}